\RequirePackage{fix-cm}
\documentclass[smallextended]{svjour3}       
\smartqed  
\usepackage{graphicx}
\usepackage{amssymb}
\usepackage{amsmath, color}
\DeclareMathOperator{\Sign}{Sign}
\def\R{\mathbb{R}}

\begin{document}

\def\a{\alpha}   \def\g{\gamma}  \def\b{\beta}
\def\l{\lambda}  \def\w{\omega}   \def\W{\Omega} \def\e{\varepsilon}
\def\d{\delta}   \def\f{\varphi}  \def\D{\Delta}    \def\r{\rho}
\def\s{\sigma}   \def\G{\Gamma}   \def\L{\Lambda}
\def\pd{\partial} \def\ex{\exists\,}  \def\all{\forall\,}
\def\dy{\dot y}   \def\dx{\dot x}   \def\t{\tau}
\def\bx{\bar x}   \def\bu{\bar u}   \def\bw{\bar w}  \def\bz{\bar z}
\def\p{\psi} \def\bh{\bar h}
\def\calO{{\cal O}}   \def\calP{{\cal P}}  \def\calD{{\cal D}}
\def\calN{{\cal N}}   \def\calL{{\cal L}}   \def\calQ{{\cal Q}}
\def\hx{\hat x}  \def\hw{\hat w}   \def\hu{\hat u}  \def\hv{\hat v}
\def\hz{\hat z}   

\def\empty{\O}
\def\lee{\;\le\;}   \def\gee{\;\ge\;}
\def\le{\leqslant}  \def\ge {\geqslant}
\def\Proof{{\bf Proof.}\,\,}

\def\ov{\overline}  \def\ld{\ldots}
\def\tl{\tilde}   \def\wt{\widetilde}  \def\wh{\widehat}
\def\lan{\langle}  \def\ran{\rangle}

\def\hsp{\hspace*{8mm}} \def\vs{\vskip 8mm}
\def\bs{\bigskip}     \def\ms{\medskip}   \def\ssk{\smallskip}
\def\q{\quad}  \def\qq{\qquad}
\def\np{\newpage}   \def\noi{\noindent}
\def\pol{\frac12\,}    \def\const{\mbox{const}\,}
\def\vad{\vadjust{\kern-3pt}}
\def\R{\mathbb{R}}   \def\ctd{\hfill $\Box$}  \def\inter{{\rm int\,}}

\def\dis{\displaystyle}
\def\lra{\longrightarrow}  \def\iff{\Longleftrightarrow}
\def\beq{\begin{equation}\label}  \def\eeq{\end{equation}}
\def\bth{\begin{theorem}\label}   \def\eth{\end{theorem}}
\def\ble{\begin{lemma}\label}     \def\ele{\end{lemma}}
\def\begar{\begin{array}} \def\endar{\end{array}}

\newcommand{\vmin}{\mathop{\rm vraimin}}
\newcommand{\vmax}{\mathop{\rm vraimax}}

\newcommand{\blue}[1]{\textcolor{blue}{#1}}
\newcommand{\red}[1]{\textcolor{red}{#1}}
\newcommand{\gr}[1]{\textcolor{green}{#1}}

\mathsurround 1pt
\def\bls{\baselineskip}   \def\nbls{\normalbaselineskip}
\bls =1.1\nbls

\title{Optimal Synthesis in a Time-Optimal Problem for the Double Integrator
    System with a Linear State Constraint }


\titlerunning{Synthesis for the Double Integrator System with a Linear State Constraint}

\author{Andrei Dmitruk        \and
        Ivan Samylovskiy 
}


\institute{Andrei Dmitruk \at
              Mathematical institute of the Russian Academy of Sciences, \\
              Russia, Moscow.\,
              \email{vraimax@mail.ru}           
        \and
        Ivan Samylovskiy \at
           Lomonosov Moscow State University, Faculty of Space Research,
           Branch of the Moscow Center for Fundamental and Applied Mathematics
           at Lomonosov Moscow State University, \\
           Russia, Moscow.\,
           \email{ivan.samylovskiy@cosmos.msu.ru}
}

\date{Received: date / Accepted: date}

\maketitle

\begin{abstract}
We consider a time-optimal problem for the classical ``double integrator''
system under an arbitrary linear state constraint. Using the maximum principle,
we construct the full synthesis of optimal trajectories and provide a qualitative
investigation of the measure, the Lagrange multiplier corresponding to the
state constraint.

\end{abstract}

\section{Statement of the Problem}
\label{sec:1}
On a time interval $[0,T]$ consider the following time-optimal problem:
\begin{equation}\label{ode}
    \begin{cases}
        \begin{aligned}
            \dot x &= y, &\q x(0) &= x_0, &\q x(T) &=0,\\[2pt]
            \dot y &=u, &\q y(0)&=y_0, &\q y(T)&=0,
        \end{aligned}
    \end{cases}
    \end{equation}
\begin{equation}\label{costu}
T\to\min, \qq |u| \le 1,
\end{equation}
under a linear state constraint
\begin{equation}\label{phase}
y\, \ge\, k x - b \qq  (b>0).
\end{equation}

Here $x$ is the position of an object (material point) on a straight line,
$y$ its velocity, the control $u$ is an acceleration applied to the object,
$t\in [0,T].$ The problem is to move the object from a given initial state $(x_0, y_0)$
to the final state $(0, 0)$ in the minimum time under the linear constraint
\eqref{phase} on the state variables $x, y.$ We assume that $u(t)$ is an
arbitrary measurable bounded function, and hence $x(t),$ $y(t)$ are Lipschitz
continuous functions.

In the absence of the state constraint, \eqref{ode}--\eqref{costu} is a well-known
Feldbaum problem that served as one of the first test examples for application of
the Pontryagin maximum principle (see \cite{Pont}). The case when the constraint
\eqref{phase} is present and $k=0,$ is considered, e.g., in the book \cite{IT}
as an example of application of the maximum principle (MP) in the form of
Dubovitskii--Milyutin.\, The case of the general constraint \eqref{phase}
has not yet been considered.

Note that \eqref{phase} is the general form of linear state constraint in this
problem which is not vertical. The constraint $y\le kx +b\;\, (b>0)$ is symmetric
to this one. The vertical linear constraints $x \le p\,$ and $x \ge -p \,\; (p>0)$
will be considered separately in Sec. \ref{sec:12} below.

If either constraint \eqref{phase} is absent, or is present with $k=0,$
the solution can be found without such an advanced theory as MP, by using just
basic tools of calculus. Indeed, here we need to find the minimal time interval
on which a Lipschitz function $y(t)$ with given endpoints conditions, the bound
on its derivative $|\dy|\le 1,$ and a lower bound on the function itself
$y \ge -b$ has a given integral. In the absence of this lower bound, we conclude
by simple considerations that the optimal function is piecewise linear with the
derivative $\pm 1$ and at most one break. If the obtained function violates the
lower bound on some time interval, then one must put $y= -b$ on this interval
and expand the latter so that the whole function $y(t)$ has the given integral.
A detailed elaboration of this solution can be recommended as an exercise for
students.

If the restriction \eqref{phase} is present\, and $k\ne0,$ the solution can hardly
be found in the described way, because here the lower bound for the function
$y(t)$ at each point $t$ depends on its integral on the interval $[0,t].$
We will apply here the maximum principle for problems with state constraints
that was obtained by A.Ya. Dubovitskii and A.A. Milyutin in \cite{DM65}.

It is well known that the ``deciphering'' of MP conditions for searching optimal
trajectories is not a routine 
task even for problems without state constraints. The more so this is true
if these constraints are present, when both the conditions and their deciphering
are more complicated. If we would knew a qualitative structure of optimal
trajectories, the search for them could be essentially more simple, but this
structure is not a priori known, and we do not impose any assumptions about that.
Therefore, we perform here a detailed procedure of solving a problem with
a state constraint by means of the MP.\,
Perhaps it can be useful for solving other problems of this type.

\section{Formulation of the Maximum Principle}
\label{sec:2}

Let a process $x^0(t), y^0(t), u^0(t),\; t\in [0, T]$ provide the minimum in
problem \eqref{ode}--\eqref{phase}. Assume that its starting point $(x_0, y_0)$
does not lie on the state boundary, i.e. $y_0 > kx_0-b.$ (The boundary case will
be considered separately.)

Then according to \cite{DM65} (see also \cite{IT,MDO,develop})\, there exist
functions of bounded variations $\varphi(t),\, \psi(t)$ (adjoint variables),
continuous from the left and {\it not equal both totally zero}, nondecreasing
function $\mu(t)$ with $\mu(0)=0,$ which generate the {\it Pontryagin function}
\begin{equation}\label{H}
H =\; \varphi y + \psi u
\end{equation}
and the {\it augmented Pontryagin function}
\begin{equation}\label{Hext}
\overline{H} =\; \varphi y + \psi u + \dot \mu \left(y - k x + b\right),
\end{equation}
so that the following conditions hold:\q the adjoint (costate) equations
\begin{equation}\label{costateEq}
\begin{cases}
\begin{aligned}
\dot \varphi & =\, -\overline{H}'_x &=&\;\; k \dot{\mu},\\[4pt]
\dot \psi & =\, -\overline{H}'_y &=&\;\, -(\varphi +\dot{\mu}),
\end{aligned}
\end{cases}
\end{equation}
the complementary slackness condition
\beq{slack}
\dot\mu(t) \left(y^0(t) - k x^0(t) + b\right)\, =\, 0,
\end{equation}
the energy conservation law
\begin{equation}
H(x^0(t), y^0(t),u^0(t))\; \equiv\; \const \ge 0,
\end{equation}
and the maximality condition
$$ \max_{|u|\le1}\,H(x^0(t),y^0(t),u)\; =\; H(x^0(t), y^0(t), u^0(t)).
$$

The latter means that
\begin{equation}\label{max}
u^0\,\in\, \Sign \psi\, := \;
\begin{cases}
\begin{aligned}
+1 ,\q & \q \mbox{  if}\q \p>0,\\[4pt]
-1 ,\q & \q \mbox{  if}\q \p<0,\\[4pt]
[-1, +1], & \q \mbox{  if}\q \p=0.
\end{aligned}
\end{cases}
\end{equation}

Here and thereafter $\dot \mu(t)$ is the derivative in the sense of generalized
functions. In other words, it is the Radon--Nikodim density of the measure
$d\mu(t)$ w.r.t. the Lebesgue measure $dt,$ so the adjoint equations
should be understood as equalities between measures, i.e., \vad
$$
d \varphi =\, k\, d{\mu}, \qq d \psi =\, -\varphi\, dt\, -d{\mu},
$$
or, in the integral form, for any $t\in [0,T]$
$$
\f(t+0) -\f(0) = \int_0^{t+0} k\,d\mu(s), \q\,
\p(t+0) -\p(0) =\, -\int_0^{t+0} \f(s)\,ds - \mu(t+0).
$$

The same relates to condition (\ref{slack}), which should be understood as the
equality $d\mu(t) \left(y^0(t) - k x^0(t) + b\right)=0.$ This means that
$d\mu(t)$ can be nonzero only on the state boundary, i.e., $d\mu(t) =0$ on any
time interval where $y^0(t) > k x^0(t) + b,$ and so, the costate variables $\p,\f$
are absolutely continuous there. \ssk

We do not write the transversality conditions, since the endpoints of
trajectory are fixed.\, The nontriviality condition means that the case
$\f\equiv \p\equiv 0$ is forbidden. In what follows, we drop the
superscript 0 for the optimal process.
\ssk

Before applying the maximum principle, we define the set of starting points
from which at least one admissible trajectory originates. This set will be denoted
by $\calD$ and called {\it admissible set}.\, Obviously, it essentially depends
on the sign of the coefficient $k.$
It follows from the Filippov theorem that for any initial position from $\calD$
there exists an optimal trajectory.

\section{The case $k>0:$ Description of the admissible set }
\label{sec:3}

Denote by $\G$ the line $y= kx-b$ (boundary of the admissible state domain),
and by $S$ the switching line in the problem without state constraint (we will
call it {\it free problem}), it is the union of two semi-parabolas:
$x= -\pol y^2,\;\; y\ge 0\,$ and $x= \pol y^2,\;\; y\le 0.$ \ssk

\begin{figure}[h]
    \begin{center}
        \includegraphics[width=0.65\linewidth,trim={0 0 140 0},clip]{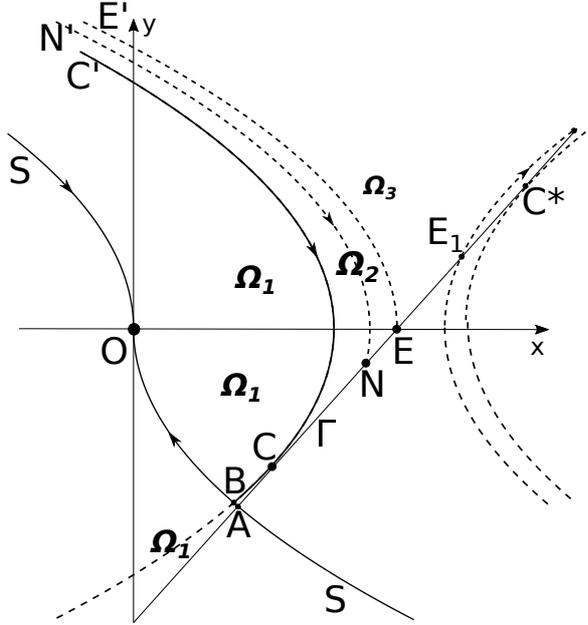}
    \end{center}
    \caption{\label{Fig1-1} Trajectory $C'CBO $ corresponds to the control
        $u= (-1,\, 1),$ trajectory $ABO$ corresponds to $u \equiv 1$}
\end{figure}

Introduce the following characteristic points for the case $k>0$\, (Fig.
\ref{Fig1-1})\footnote{We use the notation $u =(u_1, u_2)$ if
$u=u_1$ a.e. on $[0,t_1]$ and $u=u_2$ a.e. on $[t_1,T]$ for some $t_1\in (0,T).$
The similar means for the notation $u =(u_1, u_2, u_3).$}:

1) the point $A$ is the intersection of $\G$ and $S,$

2) the point $E$ is the intersection of $\G$ and the abscissa axis,

3) $C$ is the point where the line $\G$ is tangent to some parabola of the
family $x= -\pol y^2 +\const$  (clearly, this parabola is unique),

4) the point $B$ is the intersection of the tangent parabola from the previous
item with the curve $S$.
\ssk

One can easily find the coordinates of these points. The point $A$ is defined by
the relations $x= \pol y^2,\;\; y= kx-b,$ whence
$$
x_A =\, \frac{(bk +1) - \sqrt{2 bk +1}}{k^2}\,, \qq y_A=\, \frac{1- \sqrt{2 bk +1}}k\,.
$$

The point $E = (b/k,\, 0).$ The point $C$ is common to the parabola $x= -\pol y^2 + m$
and the line $y= kx-b,$ which is tangent to this parabola at $C,$ so
$dx/dy = -y = 1/k\,.\;$ Hence,
$$
y_C =\, -\frac 1 k\,, \qq x_C =\, \frac{(kb -1)}{k^2}\,,
\qq m=\, \frac{2 kb -1}{2 k^2}.
$$

The point $B$ is defined by the relations $x= -\pol y^2 + m\;\; x= \pol y^2\,,$
whence $x_B = m/2\,,$ $y_B = -\sqrt m\,.$ However, we do not need these coordinates
for qualitative research, only the relative position of the points $A$ and $C$
on the line $\G$ will be important. Depending on this relative position, there
are three possible cases, see Figs. \ref{Fig1-1}, \ref{Fig1-2} and \ref{Fig1-3},
respectively.
\ssk

First of all, consider properties of admissible trajectories that are valid
in all these three cases.

Let $x(t),\, y(t),\, u(t),\; t\in [0, T]$ be an admissible process.
For convenience, we denote by $r(t) = (x(t),\, y(t))$ the state point on
the $xOy$ plane.
We will say that a point $r= (x,y)\in\G$ lies above (below) a point
$r' =(x',y')\in\G$ and write $r > r'\,$ if $y> y'$ (respectively, $y < y').$
\ssk

1) The points on the phase boundary that lie below $C$ are never reachable.
Indeed, since $y< -1/k$ here, we have $u \ge -1> ky,$ so $(y - kx)^\bullet =$
$u - k y>0,$ hence we arrive at such a point from the forbidden area, which
is impossible.
\ssk

2) Let $C^*$ be the point on the line $\G$ that is the symmetric to $C$
with respect to $E.$ If we got to the state boundary at a point $E_1 \ge C^*,$
then further motion is impossible, because the dynamics \eqref{ode} moves us
into the forbidden area.\,
If we got to $\G$ at a point $E_1 \in (E, C^*),$ then further motion is possible
only in the lens formed by the line $\G$ and the arc of the parabola coming
from the point $E_1$ with the control $u= +1,$ and it is impossible to leave
this lens obeying the state constraint. Hence, all the points on the state
boundary above $E$ are unreachable. It follows that all the points of the upper
half-plane lying to the right of parabola $(E', E)$ and above $\G$ are also
unreachable (because from every such point we would eventually get to the state
boundary above the point $E$).
\ssk

3) An important fact is the following

\ble{lem1}
An admissible trajectory cannot pass through the point $E.$
\ele

\Proof On the contrary, suppose that $\ex t' <T$ with $r(t')=E,$ that
is $x(t')= x_E$ and $y(t')=0.$ Then, for $\bx = x- x_E$ and $\t= t-t',$ we have
$\dot{\bx} = y$ and $y\ge k\bx,$ so $\dot{\bx} = k\bx + a(\t),$ where $a(\t)\ge 0.$
Setting $\bx(\t) = c(\t)\,e^{k\t},$ we obtain $\dot c = a(\t)\,e^{-k\t}$
and $c(0)=0,$ whence $c(\t)\ge 0$ for all $\t\ge0.$ Therefore, $\bx(\t)\ge0$
and $y(\t)\ge k\bx(\t) \ge0,$ whence $\bx = y=0$ (otherwise we are in the
first quadrant and above $\G,$ which is forbidden), so in fact we stay at $E$
until the end and cannot get to the origin.
\ctd \bs

Thus, we cannot reach $E.$ But then we cannot also reach the entire upper part
of the parabola $(E',E),$ because from there we would get either to $E,$ or
to the already forbidden area on the right of this parabola.
\ssk

4) It follows that the set $\calD$ of initial positions from which there exist
admissible trajectories is the half-plane $y \ge kx-b$ except the closed set
bounded from the left by the parabola $(E', E).$
These trajectories consist of motions along the parabolas under the controls
$u= -1$ or $u= +1,$ and also motions along the state boundary.
Note that $\calD$ is not closed.
\ssk

5) The state boundary can be reached only on the semi-interval $[C,E),$ and
here $\dx = y<0.$ If at some time interval we are on $\G,$ then $y= kx-b,$
$\,\dy = k\dx <0,$ so we can move only in ``the left-down'' direction from the
point $E,$ excluding it, to the point $C.$ Since $u = ky,$ in view of the
constraint $|u|\le 1$ this motion can continue only while $|y|\le 1/k\,.$
The limit value $y = -1/k$ corresponds exactly to the point $C.$

\section{The case $k>0$:\, Analysis of MP\, }
\label{sec:4}

Let an admissible process $x(t),\, y(t),\, u(t),\; t\in [0,T]$ with a starting
point $(x_0,y_0) \notin \G$ satisfy the maximum principle.
We now pass to its analysis, which for convenience is divided into
a number of separate steps.  \ssk

\ble{lem-p}
If at a time $t'$ we are on $\G,$ then $\p(t'-0) \le 0,\,$ $\p(t'+0) \le 0,$
and $\,\D\p(t') \le 0.$

If, moreover, $r(t') >C,$ then $\D\mu(t') = 0,$ both $\p$ and $\f$ are continuous
at $t',$ and $\p(t') = 0.$
\ele

\Proof In the given case we cannot have $\p > 0$ on the left near $t',$
otherwise there would be $u= +1,$ so we came at this point from the forbidden
area. Therefore, $\p(t'-0) \le 0.$ By the costate equations,
$\D\p(t') = -\D\mu(t') \le 0,$ so $\p(t'+0) \le 0.$

If $r(t') >C$ and $\p(t'+0) <0,$ then $\p <0$ and $u= -1$ in a right neighborhood
of $t',$ so we penetrate the boundary.\, Then $\p(t'+0)\ge 0$ and $\D\p(t')\ge 0.$
But since we already proved $\D\p(t') \le 0,$ we get $\D\p(t') = 0,$ so the
measure has no jump, and $\p(t') = 0.$ \ctd
\ssk

Thus, a jump of measure is possible only at the point $C.$
\ssk

\ble{lem-pwf}
If $\p=0$ on some open interval $\w = (t',t''),$ then $\f$ is a nonpositive
exponential function there.
\ele

\Proof\, In this case $\dot\p = -(\f+\dot\mu) =0,$ whence $\f= -\dot\mu \le 0$
and $\dot\f= k\dot\mu= -k\,\f,$ hence $\f$ is an exponential function on $\w$. \ctd
\ssk

\ble{mint}
If on some open interval $(t_1, t_2)$ we are on $\G,$ then $\p=0$ and the
measure has a density $\dot\mu \ge 0$ which is an exponential function.
Moreover, $\D\mu(t_1)= 0,$ both $\p$ and $\f$ are continuous at $t_1\,,$
so they may have jumps only at $t_2$.
\ele

\Proof\, According to Step 5, $u = ky <0$ a.e. on $(t_1\,,t_2),$ i.e. we move
in the ``left-down'' direction along the interval $(E,C),$ hence $r(t) >C,$
$|y| <1/k\,$ for $t<t_2\,,$ so $|u|<1$ and by \eqref{max} $\p =0.$
Then by Lemma \ref{lem-p} $\D\mu(t_1)= 0,$ whence by the costate equations
$\p,\, \f$ are continuous at $t_1\,,$ and by Lemma \ref{lem-pwf}
$\dot\mu = -\f$ is an exponential function. \ctd
\ssk

\ble{impos}
The following case is impossible:\, $t_1 <t_2\,,$ $r(t_2)\in\G,$
we are outside of $\G$ for $t <t_1$ and $t >t_2\,,$ both $\p$ and $\f$ are
continuous at $t_1$ and vanish on $(t_1, t_2).$
\ele

\Proof\, In this case $\f(t_1) = \p(t_1) =0$ with $\dot\f= 0$ and $\dot\p = -\f$
on $[0,t_1],$ so we obtain $\f = \p= 0$ there, hence also on $[0,t_2).$
Let us now see what is on the right of $t_2\,.$ If $\D\mu(t_2)= 0,$
then $\p = \f = 0$ there, hence on the whole $[0,T],$ which is the trivial
collection of multipliers. Then $\D\mu(t_2)> 0,$ so $r(t_2) =C,$
moreover $\D\f(t_2)> 0,$ $\D\p(t_2)< 0,$ and $\dot\p = -\f <0$ for $t>t_2\,,$
so $\p<0$ and $u= -1$ for all $t>t_2\,.$
But we cannot get to the origin with this control. \ctd
\ms

6) Let $M = \{\,t\;|\; r(t)\in\G\}$ be the set of contact points with the state
boundary. If it is empty, the solution is well known --- it is the motion along
the parabolas of the two above families, see \cite{Pont}.\,

Assume that $M$ is nonempty.\, Denote $\,t_1 =\, \min M,$ $\; t_2 =\, \max M.$

\ble{lem2}
$M$ is connected, i.e. is a closed interval or a point: $M = [t_1, t_2].$
\ele

\Proof Suppose on the contrary that $M$ is not connected, i.e. there exist points
$t'< t''$ in $M,$ such that on $\w = (t',t'')$ we are outside of $\G.$
Then $\dot\mu =0$ on $\w,$ hence $\f = \const,$ and $\p$ is a linear function.\,
By Lemma \ref{lem-p}, $\p(t'-0) \le 0$ and $\p(t''+0) \le 0.$
If $\p<0$ on $\w,$ then $u= -1,$ and we will not get back to $\G$ at $t''.\,$
So, $\p =0$ on $\w.$

It follows that on any open interval $\w'\subset [t_1, t_2]$ where we are
outside of $\G,$ we have $\p =0$ (since $\w'$ is contained in some maximal
interval $\w$ of the above type), and then, due to the control system, $\f=0$
as well. Thus, if on some interval in $[t_1, t_2]$ we have $\p<0,$ then
we are on $\G$ on this interval.
\ssk

We claim that $\p =\f= 0$ on $(t_1, t_2).$ Indeed, suppose $\p(t_*) < 0$
at some point $t_*\in (t_1, t_2)$ of its continuity (the case $\p(t_*) > 0$
is excluded by Lemma \ref{lem-p}). Then $\p<0$ in a neighborhood of this point,
so we are on $\G$ and $u= -1.$ But it is impossible to keep on $\G$ with this
control. Therefore, $\p(t) =0$ at all points of its continuity in $(t_1, t_2),$
and hence, everywhere on $(t_1, t_2).\,$ Then by Lemma \ref{lem-pwf} $\f$
is an exponential function. Since $\f=0$ on the above interval $\w,$
we get $\f=0$ on $(t_1, t_2)$ as well.

Since $\p(t_1+0)=0,$ while by Lemma \ref{lem-p} $\p(t_1-0)\le 0$ and
$\D\p(t_1) \le 0,$ we get $\p(t_1-0) = 0,$ so $\D\p(t_1)= 0$ and the functions
$\p,\,\f$ are continuous at $t_1$ with $\p(t_1) = \f(t_1)= 0.$
We came to the situation of Lemma \ref{impos}, which is impossible. \ssk

It follows from these arguments that an interval $\w = (t',t'')$ with the
above properties cannot exist, and so $M = [t_1, t_2].$ The lemma is proved. \ctd
\ms

Now, consider all possible cases of location of the set $M.$ This consideration
depends on the relative positions of the points $A$ and $C.$
\ssk

\subsection{The case $k>0$ and $C> A$ }
\label{sec:5}

\noi
Define the following closed sets (see Fig. \ref{Fig1-1}):
\ssk

$\W_1$ lying on the left of the line $\G$ and the parabola $(C',C],$

$\W_2$ bounded by the parabolas $(C',C],\; (E',E]$ and the line segment $[C,E],$

and $\W_3$ lying between the parabola $[E,E')$ and the line $\G$.
\ssk

\noi
Recall that domain $\W_3$ is forbidden, so the set of admissible starting
positions $\calD =\, (\W_1 \cup \W_2) \setminus \W_3$.
\ssk

7) Can we have $M= \{t_*\}$ (touching the boundary at a single point) with
$r(t_*) >C$?\, By Lemma \ref{lem-p}, there is no jump here, so $\dot\mu \equiv 0,$
and the motion is the same as in the free problem, i.e. along the parabolas.
But this motion has $u= -1$ in a neighborhood of the point $r(t_*),$ then
moving from this point we violate the boundary. So, this case is not possible.
\ssk

8) Can we have $M = [t_1, t_2]$ with $t_1< t_2$ and $r(t_2) >C$?\,
By Lemma \ref{mint}, in this case there is no jump anywhere, $\p=0$ and
$\dot\mu = -\f\ge 0$ is an exponential function on $M.$
The case $\f=0$ is excluded by Lemma \ref{impos}, so $\f<0$ on $M.$

For $t> t_2$ the measure switches off, $\f =\f(t_2) <0,\,$ $\dot\p = -\f >0,$
then $\p >0,$ $u= +1$ until the end, and we cannot get to the origin with this
control from $r(t_2) >C.$\, So, this case is impossible.

\ssk

9) Suppose that $M = \{t_*\}$ and $r(t_*)=C.$ The first thing that comes to
mind in this case is the trajectory $C'CBO$ starting at a point $C'$ on the
parabola $C'C$ and having the control $u= (-1, -1, +1).$
Anyway, let us find what the MP implies.

Before and after the point $t_*\,,$ the function $\p$ is linear (in general,
with different coefficients), and by Lemma \ref{lem-p} $\p(t_*-0)\le 0$ and
$\,\p(t_*+0)\le 0.$\,
If $\p(t_*+0)=0,$ then $\D\mu(t_*)=0\,$ and $\p$ is totally linear
with $\p(t_*)=0.$ The case $\p\equiv 0$ implies also $\f \equiv  0,$
which is forbidden. Then, either $\p<0,$ $u= -1$ for all $t> t_*\,,$
or $\p>0,$ $u= +1$ there. In both of these cases we do not reach the origin.

Thus, $\p(t_*+0)< 0$ and $\p$ increases on $[t_*, T]$ (so $\f<0$ there) changing
its sign from minus to plus at a time $t_B \in (t_*, T)$ corresponding to the
point $B,$ with the respective switching of the control from $-1$ to $+1.$\,

Since $\dot\mu=0$ both on the left and right of $t_*\,,$ then
$\D\dot\p(t_*)= -\D\f(t_*) \le 0,$ so $\p$ increases as well on $[0,t_*)$
up to a value $\p(t_*-0)\le 0,$ hence $\p<0$ and $u= -1$ on this interval.
We thus move along the above trajectory $C'CBO.$
Since it is optimal in the free problem, it just occasionally touches
the state boundary, which makes no affection on it.

However, let us find a possible value of the jump $\D\mu(t_*) := \d$ for this
trajectory. Denote by $t_C\,,\, t_B$ the moments corresponding to the points
$C, B$ (here $t_C =t_*\,).$
Choose the normalization $\dot\p = -\f= 1$ on $(t_C,T]$ for our collection
of multipliers. Then $\p =\, t- t_B$ on $(t_C, T]$ independently of $\d,$ and
we should have $\p(t_C-0)= (t_C- t_B) +\d\le 0,$ whence $0\le\d\le t_B -t_C\,.$
For any such $\d$ we have $\dot\p(t_C-0) = -\f(t_C-0) = 1 + k\d > 0,$
hence $\p<0$ on $[0,t_C],$ so the MP is satisfied.\,
Thus, a possible value of the jump of measure at the point $C$ {\it is
not unique}, 
being arbitrary within the above bounds\footnote{As a test, note that
for any jump of measure we have
$\D H(t_*) =\D\f\cdot y(t_*)+|\D\p| =k\,\d\,(-\frac 1 k)\,+\,\d =0,$
so the condition $H= \const$ is valid independently on the value of $\d.$}.
In particular, if $\d=0,$ then $\p$ is totally linear. If $\d =1,$ then $p<0$
on $[0,t_C]$ and linearly increases up to the value $\p(t_C-0) =0.$

Thus, the case when $M = \{t_*\}$ and $r(t_*)=C>A$ is completely done.
\ssk

10) Finally, let $M = [t_1, t_2],$  $t_1< t_2$ and $r(t_2) =C.$
Then, by Lemma \ref{mint}, within $M$ we move along the boundary in the
``left-down'' direction, the jump at $t_1$ is absent, $\p$ is continuous at
this point and vanishes on $[t_1, t_2),$ while $\f = -\dot\mu \le 0\,$ is
an exponential function there. \ssk

The case when $\f=0$ on $M$ is excluded by Lemma \ref{impos}, so $\f <0$
on $M.$ Moreover, this holds on $[0,t_1]$ as well. Then $\dot\p = -\f >0,$
so $\p <0$ and increases up to zero on $[0,t_1],$ while $u= -1.$

If the jump at $t_2$ is absent, the measure switches off on $[t_2,T],$ then
$\f = \const <0,\,$ $\p(t_2) =0$ and $\dot\p = -\f >0,$ whence $\p >0$ and
$u= +1$ until the end of motion.  Since $r(t_2)=C >B,$ we cannot get
to the origin under this control.

Thus, the measure must have a jump $\D\mu(t_2) = \d >0.$
Then $\D\f(t_2) = k\d$ and $\p(t_2+0)= -\d.$
For definiteness, choose the normalization $\dot\p= -\f= 1$ for $t> t_2\,.$
Then $\f(t_2-0) = -1 -k\d,$ which gives the above relations $\dot\p>0,$ $\p<0,$
and $u= -1$ on $[0,t_1].$ At the switching point $B$ we should have
$\p(t_B)= -\d + \D t =0,$ where $\D t = t_B -t_C = y_C -y_B$ is the motion time
along the arc $[CB]$ with $u = -1.$ Therefore,  $\d =\, \D t.$
For the given positions of the points $C$ and $B$ (which are determined by the
value of $k),$ the value of $\D t$ is known, so the value of jump $\d$
is {\it uniquely determined}.\, So, we have a trajectory of the type $N'NCBO$
with $u = (-1,\,ky,\, -1, +1)$ that starts at a point $N' \in \inter\W_2\,,$
the first leg of which is the motion with $u=-1$ to a point $N \in (C,E).$
\ssk

11) Thus, for the case $k>0$ and $C > A,$ we examined all possible cases of
location of the set $M$ (the contact set with the state boundary), and for
each case we specified the corresponding costate variables and the measure
supported on this set. \ssk

Since for any starting point in $\calD$ there is only one trajectory satisfying
the MP, it is indeed optimal.\, This yields the following optimal synthesis.
\ssk

\bth{synth1}
Let $k>0$ and $C>A$ (Fig. \ref{Fig1-1}).
In the domain $\W_1$ the optimal trajectories are the same as in the free
problem.

In the domain $\inter\W_2$  any optimal trajectory is of the type $N'NCBO$
with the control $u= (-1, ky, -1, +1)$ and a unique positive jump of measure
at the point $B,$ where $N$ is a point on the line interval $(C,E).$

For starting points lying on the interval $(C,E),$ the optimal trajectories
are the ``tails'' of optimal trajectories starting from $\inter\W_2\,,$
the first leg of which is the motion with $u= -1$ to a point on $(C,E).$
\eth

\subsection{The case $k>0$ and $C =A=B$ }
\label{sec:5}

\begin{figure}[h]
    \begin{center}
        \includegraphics[width=0.7\linewidth,trim={0 0 0 0},clip]{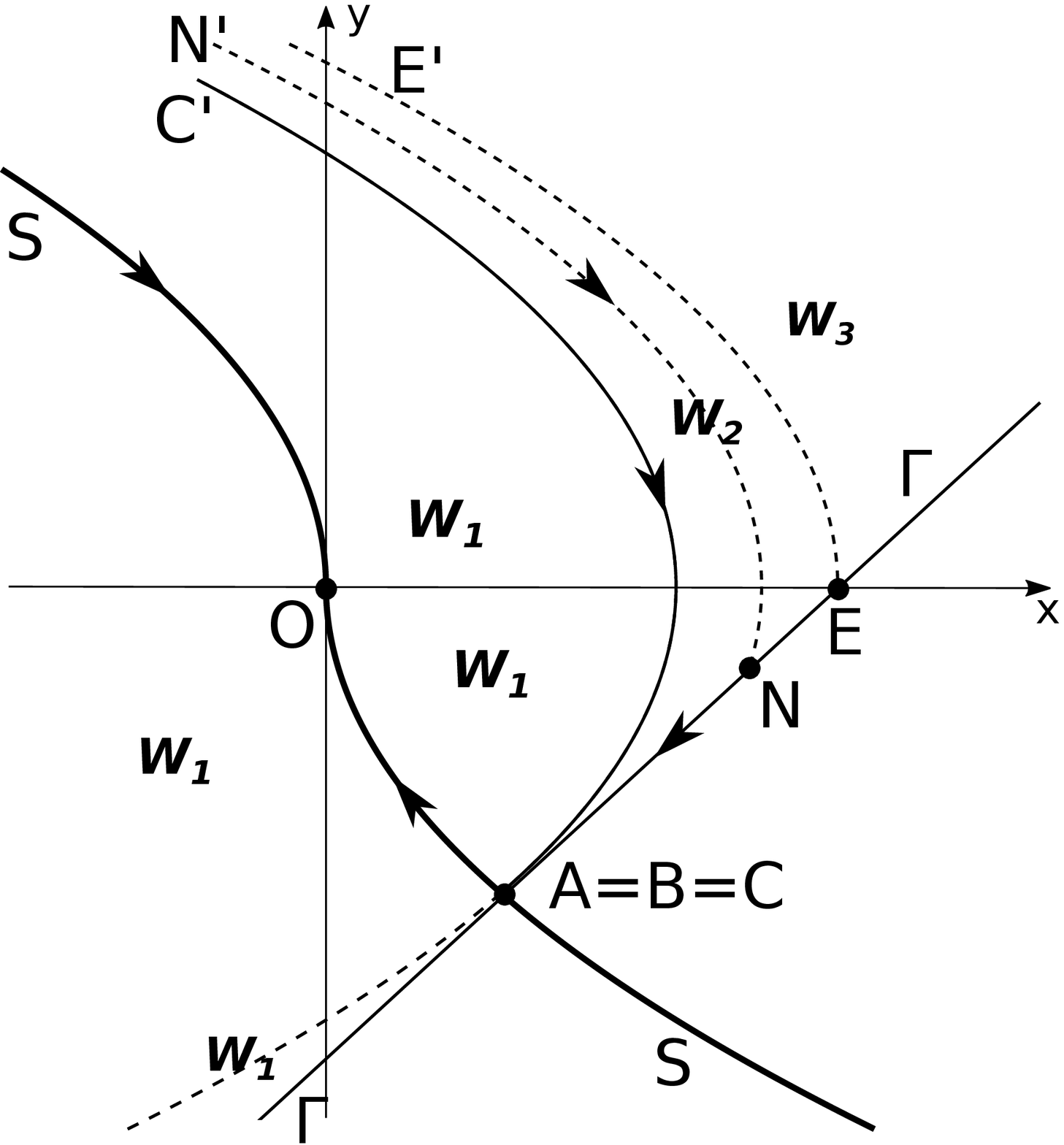}
    \end{center}
    \caption{\label{Fig1-2}}
\end{figure}

\noi
Define the following closed sets (see Fig. \ref{Fig1-2}):
\ssk

$W_1$ lying on the left of the line $\G$ and the parabola $(C',C],$

$W_2$ bounded by the parabolas $(C',C],\; (E',E]$ and the line segment $[C,E],$

and $W_3$ lying between the parabola $[E,E')$ and the line $\G$ (forbidden).
\ssk

\noi
Here, the admissible set $\calD = (W_1 \cup W_2) \setminus W_3\,.$
One can see that Steps 7 and 8 are still valid (see Fig. \ref{Fig1-2}).
\ms

12) Consider the case when $M = \{t_*\}$ and $r(t_*)=C.$ Like in Step 9,
before and after $t_*$ the function $\p$ is linear (may be, with different
coefficients), $\p(t_*-0)\le 0,$ $\,\p(t_*+0)\le 0,\,$ $\D\p(t_*)\le0,\,$
and $\D\dot\p(t_*) \le 0.$ \ssk

If $\D\mu(t_*)>0,$ then $\p<0$ on the right near $t_*\,,$ and further
may change its sign to plus, so either $u= -1$ or $u= (-1,+1)$ on $[t_*,T].$
However, in both these cases we do not reach the origin. Therefore, $\D\mu(t_*)=0,$
whence $\p$ is totally linear and obviously increases, so $u= (-1,+1)$ on $[0,T]$
with the switching point $t_*\,.$ We obtain the trajectory $C'CO$ which is
optimal in the free problem. It just occasionally touches the state boundary,
which makes no affection both on the trajectory and on the costate variables.

13) Now, let $M = [t_1, t_2],\,$ $t_1< t_2$ and $r(t_2) =C.$
Like in Step 10, within $M$ we move along the boundary in the ``left-down''
direction, $\p <0$ on $[0,t_1]$ and increases up to zero, while $u= -1.$
Further on, $\p$ is continuous at $t_1$ and vanishes on $[t_1, t_2);\,$
moreover, $\f <0$ on $[0,t_2).$

If $\D\mu(t_2)>0,$ then $\p<0$ on the right near $t_2$ and further may
change its sign to plus, so either $u= -1$ or $u= (-1,+1)$ on $[t_2,T].$
But in both these cases we do not reach the origin. Therefore, $\D\mu(t_2)=0$
(the measure is absolutely continuous), hence for $t> t_2$ we have $\f <0,$
so $\p>0,$ $u= +1,$ and we obtain the trajectory $N'NCO$
with $u = (-1,\, ky,\, +1).$
\ssk

Thus, for any starting point in $\calD,$ a trajectory satisfying the MP
is unique, so it is optimal.\, Therefore, we obtain

\bth{synth2}
Let $k>0$ and $C=A=B$ (Fig. \ref{Fig1-2}). In the domain $W_1$ the optimal
trajectories are the same as in the free problem.

In the domain $W_2$ an optimal trajectory is of the above type $N'NCO$
with $u = (-1,\, ky,\, +1),$ where $N$ is a point on the line interval $(E,C).$

For starting points lying on the interval $(C,E),$ the optimal trajectories
are the "tails"\, of optimal trajectories starting from $\inter W_2\,,$
the first leg of which is the motion with $u= -1$ to a point on $(C,E).$
\eth

\subsection{The case $k>0\,$ and $C < A$}
\label{sec:6}

\begin{figure}[h]
    \begin{center}
        \includegraphics[width=0.6\linewidth,trim={0 0 0 0},clip]{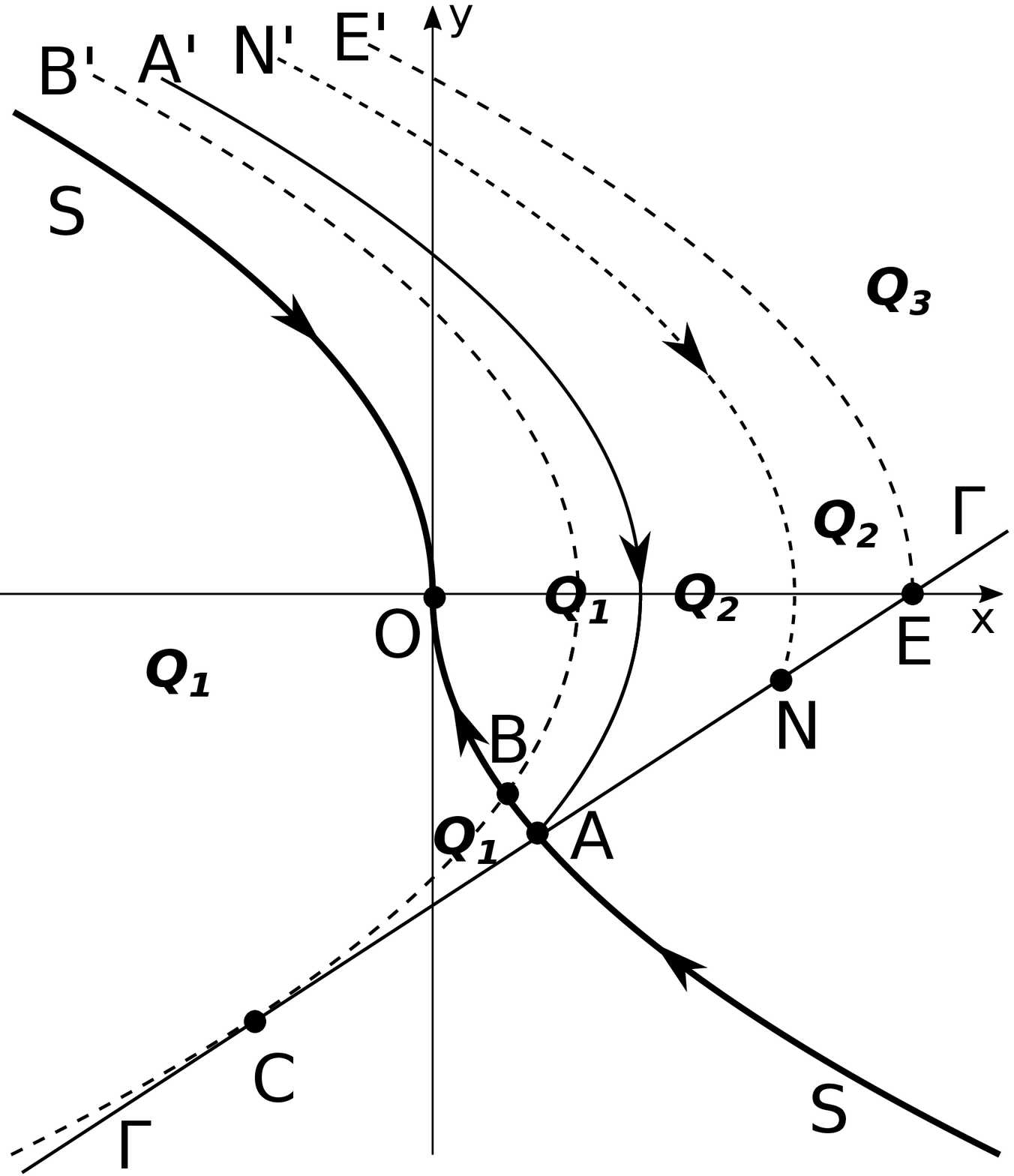}
    \end{center}
    \caption{\label{Fig1-3}}
\end{figure}

\noi
Define the following closed sets (see Fig. \ref{Fig1-3}):
\ssk

$Q_1$ lying on the left of the line $\G$ and the parabola $(A',A],$

$Q_2$ bounded by the parabolas $(A',A],\; (E',E]$ and the line segment $[A,E],$

and $Q_3$ lying between the parabola $[E,E')$ and the line $\G$ (forbidden).
\ssk

\noi
The admissible set $\calD = (Q_1 \cup Q_2) \setminus Q_3\,.$
This case differs from the case $C\ge A$ in that here the motion along the
boundary $\G$ below the point $A$ is of no sense, because one can switch to
$u =+1$ at this point and get to the origin in a smaller time.

Therefore, by Lemma \ref{lem-p} the measure has no jumps here altogether.
By analogy to Step 7, one can see that the case when $M = \{t_*\}$ with
$r(t_2) >A$ is impossible, and by analogy to Step 8, the case when $M =[t_1,t_2]$
with $t_1 <t_2$ and $r(t_2) >A,$ is also impossible.\,

If $M = \{t_*\}$ with $r(t_*) =A,$ the measure is absent by Lemma \ref{lem-p},
so we have the trajectory $A'AO$ of the free problem.
\ssk

14) Consider the only possible case with a nonzero measure, which is when
$M = [t_1, t_2]$ with $t_1 <t_2$ and $r(t_2) =A.$ Then, by Lemma \ref{mint}
$\p=0$ on $M,$ while by Lemma \ref{impos} $\f<0$ on $M$ and hence on the whole
interval $[0,T].$ Since $\dot\p= -\f>0$ both on the left of $t_1$ and the right
of $t_2\,,$ we have $\p<0,\;\, u= -1$ and $\p>0,\;\, u= +1$ on these intervals,
respectively.\, Therefore, we obtain a trajectory $N'NAO$ with $u= (-1,\,ky,\,+1),$
where $N$ is a point on the line interval $(E,A).$
\ssk

Since for any starting point in $\calD,$ a trajectory satisfying the MP
is unique, it is optimal.\,
Thus, here the optimal synthesis is as follows.

 \ms

\bth{synth3}
Let $k>0$ and $C < A$ (Fig. \ref{Fig1-3}).
In the domain $Q_1$ the optimal trajectories are the same as in the free problem.

In the domain $Q_2$ any optimal trajectory is of the above type $N'NAO$
with $u = (-1,\, ky,\, +1),$ where $N$ is a point on the line onterval $(A,E).$

For starting points lying on $(A,E),$ the optimal trajectories are the ``tails''
of optimal trajectories starting from $\inter Q_2\,,$
the first leg of which is the motion with $u= -1$ to a point on $(A,E).$
\eth

{\it Remark 1.}\, Note an important property of the Bellman function $V(x_0, y_0)$
for the studied problem, that appears in both three subcases for $k>0.$
Recall that it is the minimal time of motion when starting from a point
$(x_0, y_0).$ Obviously, $V(x_0, y_0)$ is finite for any point in $\calD.$
However, if the point $(x_0, y_0)$ tends to boundary parabola $(E',E],$ then
$V(x_0, y_0)\to +\infty.$ This is caused by the fact that if $(x_0, y_0)$ is
close to that parabola, then at some moment $t_*$ the corresponding optimal
trajectory gets to a point $r(t_*)$ of the interval $(C,E)$ and then goes along
this interval up to the point $C$ or $A.$ Since the latter motion obeys
$\dy = k y,$ the time of motion along the interval $[r(t_*),C]$ or $[r(t_*),A]$
tends to infinity as $r(t_*)$ tends to $E.$

\section{The case $k<0$}
\label{sec:7}

Denote $k = -q,$ where $q>0.$ The boundary $\G$ of the admissible state
domain is now defined by $y= -qx-b,$ while the switching line $S$ is the same
as before. There are three qualitative different cases here, depending on
the number of points where the line $\G$ and the curve $S$ intersect:\,
one, two, or three. 

\subsection{The case $k<0$ with three intersection points}
\label{sec:8}

\begin{figure}[h]
    \begin{center}
        \includegraphics[width=0.8\linewidth,trim={0 0 0 0},clip]{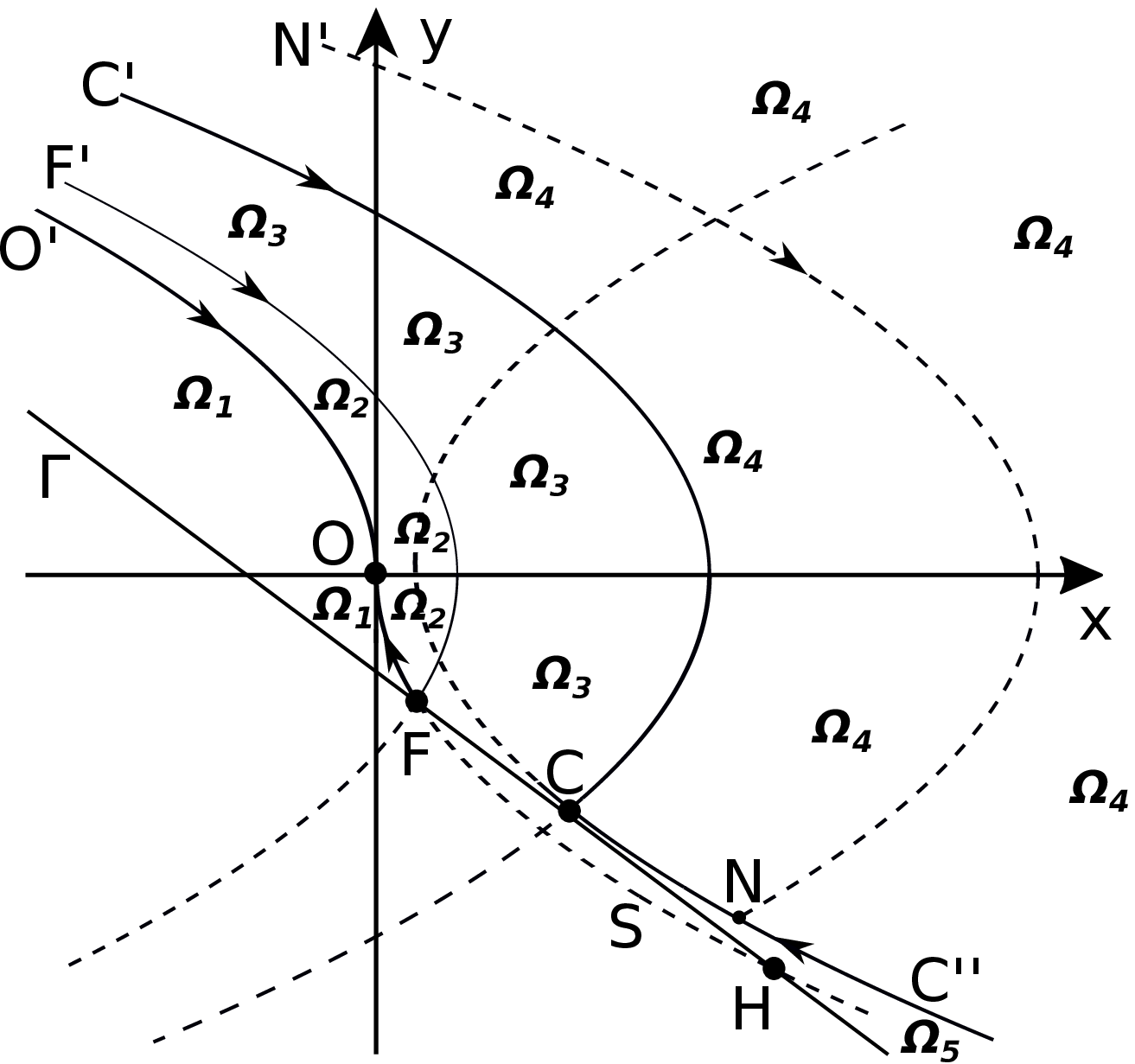}
    \end{center}
    \caption{\label{Fig2}}
\end{figure}

Introduce characteristic points for this case. Let $C$ be the point where
the line $\G$ is tangent to some parabola of the family $x= \pol y^2 +\const$
(clearly, this parabola is unique), $F$ be the intersection point of $\G$
and $S$ lying above the point $C,$ and  $H$ be the intersection point of $\G$
and curve $S$ lying below the point $C$ (see. Fig. \ref{Fig2}). \ssk

Let us find the coordinates of these points. The point $F$ is defined by the
relations $\,x= \pol y^2,$ $\; y= -(q x+ b),$ whence
$$
x_F =\, \frac{(1- bq) -\sqrt{1-2 bq}}{q^2}\,, \qq y_F =\, \frac{-1+ \sqrt{1-2 bq}}q\,.
$$
The point $C$ is common to a parabola $x= \pol y^2 + m$ and the line
$y= -(qx +b),$ which is tangent to this parabola at $C,$ so $dx/dy = y = -1/q\,.$
Hence,
$$ y_C =\, -\frac 1 q\,, \qq x_C=\, \frac{(1- bq)}{q^2}\,. 
$$
The point $H$ will not play any role in our analysis, as well as the point
of intersection of $\G$ with the parabola $[O,O')$ lying in the second quadrant.
\ssk

\noi
Define the following closed sets:

$\W_1$ bounded by the line $\G$ and the left branch of line $S,$

$\W_2$ lying between the left branch of line $S$ and the parabola $[F,F'),$

$\W_3$ bounded by the line $\G$ and the parabolas $[F,F')$ and $[C,C'),$

$\W_4$ bounded by the parabolas $[C,C')$ and $[C,C''),$

and $\W_5$ lying between the line $\G$ and parabola $[C,C'').$
\ssk

\noi
First of all, consider some properties of admissible trajectories of the control
system \eqref{ode} valid in our case, and determine the set of admissible
starting positions.
\ssk

1) The points on the phase boundary lying below $C$ are not reachable on
the trajectories of our system.\, Indeed, here always $y< -1/q,$ so
$(y+ qx)^\bullet = u+ q y< u-1\le 0,$ i.e. the admissible velocities are
directed inside the forbidden area, whence further movement from such points
is impossible.

It is easily seen that any points of $\inter\W_5$ are also inadmissible,
because any motion from them gets us to the phase boundary below $C,$
so the set $\calD$ of admissible starting positions is the one bounded from
below by the line $\G$ and parabola $[C,C'').$ In other words,
$\calD = \W_1 \cup \W_2 \cup \W_3 \cup \W_4\,,$ and it is a closed set.

\ssk

2) Observe that for the initial points lying in the regions $\W_1$ and $\W_2$
the optimal trajectories of the free problem, without the state constraint,
satisfy this constraint, so they are optimal in the problem with this constraint
as well. Thus, it remains to consider only the regions $\W_3$ and $\W_4\,.$
\ssk

3) According to Step 1, the state boundary can be reached only at the points
above or equal the point $C,$ where $\dx = y<0.$
If on some time interval we are on the boundary $\G,$ i.e. $y(t)= -qx(t)-b,$
then $\dy(t) = -q\,\dx(t) >0,$ so we move in ``the left-up'' direction,
as long as the condition $y<0$ holds. Clearly, it is enough to move only to
the point $F,$ because from it we can optimally get to the origin with
the control $u= +1.$ Since $u = -q\,y,$ in view of the constraint $|u|\le 1$
this motion can happen only while $y \ge -1/q\,.$ The limit value $y = -1/q $
corresponds to the point $C.$ Thus, the motion along the state boundary
is possible only within the segment $[C, F].$

\subsection{Analysis of MP for the case $k<0$ with three intersection points}
\label{sec:9}

Let $(x(t),\, y(t),\, u(t)),\;\, t\in [0,T],$ be an optimal process starting
at a point $(x(0), y(0)) \in \W_3 \cup \W_4$ not lying on $\G$
(see Fig. \ref{Fig2}). Then it satisfies the maximum principle
\eqref{costateEq}--\eqref{max}. We now pass to its analysis.
Partially, it will be similar to that in the case $k>0.$

\ble{sign-p}
If at a time $t'$ we are on $\G,$ then $\p(t'-0) \ge 0,\,$ $\p(t'+0) \ge 0,\,$
$\D\p(t') \le 0,$ and $\D\dot\p(t') \ge 0.$

If, moreover, $r(t') >C,$ then $\D\mu(t') = 0,$ both $\p$ and $\f$ are continuous
at $t',$ and $\p(t') = 0.$
\ele

\Proof\, In a right neighborhood of $t'$ we cannot have $\p < 0,$
otherwise there would be $u= -1,$ and we violate the state boundary.
Then $\p(t'+0) \ge 0.$ By the costate equation $\D\p(t')= -\D\mu(t')\le 0,$
so $\p(t'-0) \ge 0,$ and also $\D\dot\p(t')= -\D\f(t')\ge 0,$  \ctd

If $r(t')> C$ and $\p(t'-0) >0,$ then $p>0$ and $u= +1$ in a left neighborhood
of $t',$ so we came at this point from the forbidden area.
Therefore, $\p(t'-0) \le 0$ and $\D\p(t')\ge 0.$ But we already have
$\D\p(t')\le 0,$ which yields $\D\p(t')= 0,$
the measure has no jump, and $\p(t') = 0.$ \ctd
\ssk

Thus, a jump of measure is possible only at the point $C.$
\ssk

Lemma \ref{lem-pwf} is still valid in this case, while Lemmas \ref{mint}
and \ref{impos} transform into the following ones.

\ble{mint2}
If on some open interval $(t_1, t_2)$ we are on $\G,$ then $\p=0$ and the
measure has a density $\dot\mu \ge 0$ which is an exponential function.
Moreover, $\D\mu(t_2)= 0,$ both $\p$ and $\f$ are continuous at $t_2\,,$
so they may have jumps only at $t_1$.
\ele

\Proof\, According to Step 3, $u = -q y >0$ a.e. on $(t_1\,,t_2),$ i.e. we move
in the ``left-up'' direction along the interval $(C,F),$ hence $r(t) >C,$
$|y| <1/k\,$ for $t >t_1\,,$ so $|u|<1$ and by \eqref{max} $\p =0.$
Then by Lemma \ref{sign-p} $\D\mu(t_2)= 0,$ whence by the costate equations
$\p,\, \f$ are continuous at $t_2\,,$ and by Lemma \ref{lem-pwf}
$\dot\mu = -\f\,$ is an exponential function. \ctd
\ssk

\ble{impos2}
If the initial point does not lie on the parabola $(C,C''),$
the following case is impossible:\, $t_1 <t_2\,,$ $r(t_1)\in\G,$
we are outside of $\G$ for $t <t_1$ and $t >t_2\,,$ both $\p$ and $\f$ are
continuous at $t_2$ and vanish on $(t_1, t_2).$
\ele

\Proof\, In this case $\f(t_2) = \p(t_2) =0$ with $\dot\f= 0$ and $\dot\p = -\f$
on $[t_2,T],$ so we obtain $\f = \p= 0$ there, hence also on $[t_1, T).$
Let us now see what is on the left of $t_1\,.$ If $\D\mu(t_1)= 0,$
then $\p = \f = 0$ there, hence on the whole $[0,T],$ which is the trivial
collection of multipliers. Therefore, $\D\mu(t_1)> 0,$ whence $\D\f(t_1)< 0,$
$\D\p(t_2)< 0,$ and $\dot\p = -\f <0$ for $t>t_2\,,$ so $\p>0$ and $u= +1$
for all $t< t_1\,.$  By Lemma \ref{sign-p}, $r(t_1) =C$ so the initial
point lies on the parabola $(C,C'')$ on the right of $C,$ which is excluded
from consideration. \ctd
\ms

6) Consider the set $M$ of contact points with the state boundary.
If it is empty, we have the motion along the parabolas of the above two
families, which is possible only for the regions $\W_1$ and $\W_2\,.$\,
We then assume that $M$ is nonempty. Denote $t_1 =\min\,M,$ $\; t_2 =\max\,M.\,$
The following lemma is similar to Lemma \ref{lem2}, but refers to
another situation, depicted in Fig. \ref{Fig2}.

\ble{lem3}
If the initial point does not lie on the parabola $(C,C'')$ on the right of $C,$
the set $M$ is connected, i.e. it is a closed interval or a point: $M = [t_1, t_2].$
\ele

The proof in fact repeats that of Lemma \ref{lem2} with corresponding
alterations:\, the points $t_1$ and $t_2$ change their roles, so do the
intervals $[0,t_1]$ and $[t_2,T].$ We live the details to the interested reader.
\ssk

Like before, consider all possible cases of location of the set $M.$
\ssk

8) Can we have $M= \{t_*\}$ (touching the boundary at a single point) with
$r(t_*)<F$?\, If $\D\mu(t_*)= 0,$ then $\dot\mu \equiv 0,$ so the
motion is the same as in the free problem, i.e. along the parabolas.
But since $r(t_*)< F,$ this motion has $u= -1$ in a neighborhood of $r(t_*),$
and we penetrate the boundary. Therefore, $\D\mu(t_*)>0,$ which is possible
only if $r(t_*) =C.$

In this case, the function $\p$ is linear on the right of $t_*$ (and also on
the left, with another slope), and since by Lemma \ref{sign-p} $\p(t_*+0) \ge 0,$
we get either $\p>0$ on the right of $t_*\,$ and then $u= +1,$ or $\p$
changes its sign there from plus to minus and then $u= (+1, -1).$
It is easily seen that in both of these cases we cannot get to the origin from
the point $C,$ so they are impossible.  \ssk

Thus, the set $M$ can consist of a single point only if this point is $F,$ and
by Lemma \ref{sign-p}\, the measure has no jump here. Then we have a trajectory
of the free problem just occasionally touching the state boundary.
\ssk

9) Can we have $M = [t_1, t_2]$ with $t_1 <t_2$ and $r(t_1)>C$? \,
In this case there are no jumps of the measure anywhere, and by Lemma \ref{mint2}\,
$\p =0$ on $M,$ so $\f = -\dot\mu \le 0$ is an exponential function.
If $\f=0$ on $M,$ then $\f =\p =0$ everywhere, which is forbidden.
On the right of $t_2$ we have $\f =\const <0, $ then $\dot\p= - \f> 0,$ hence
$\p> 0, $ $\,u = 1$ without switching, and then we can get to the origin only
if $ r(t_2) = F.$ Thus, the boundary arc of the trajectory always ends at $F.$

On the first segment $[0, t_1]$ we also have $\f =\const <0 $ and
$\dot\p = -\f> 0,$ hence $\p <0$ and grows up to the value $ \p(t_1)= 0$
with $u = -1.$ Such trajectories correspond to the starting positions
from the region $\W_3 \,. $
\ssk

10) Now, let $M = [t_1, t_2],$ where $t_1< t_2\,,$ $r(t_1) =C,$
and $r(t_2) = F.\,$ Assume also that the starting position $r(0)$ does not lie
on parabola $(C, C '').\,$
By Lemma \ref{mint2}, $\p=0$ on $(t_1, t_2]$ and $\f= -\dot\mu$ is an exponential
function. The case when also $\f=0$ on $(t_1, t_2]$ is excluded by Lemma \ref{impos2}.
Then $\f<0$ on $(t_1, t_2]$ and thus on $(t_1, T].$

For $t> t_2$ the measure switches off and we have $\dot\p = -\f >0,$ so
$\p >0,$ and we move with $u= +1$ along the arc $FO$ until the end.
\ssk

If $\D\mu(t_1) =0,$ then $\f =\const <0$ on the left of $t_1,$ and
$\dot\p = -\f >0,$ hence $\p<0$ and increases to the value $\p(t_1)=0$
while $u = -1.$ This corresponds to the trajectories starting on parabola $(C,C')$
that separates regions $\W_3$ and $\W_4\,.$ \ssk

If $\D\mu(t_1)>0,$ then $\D\p(t_1)< 0$ and $\p(t_1-0) >0.$ In this case,
either $\p >0$ everywhere on $(0,t_1],$ or $\p$ changes its sign from minus
to plus. The first case corresponds to the trajectories starting from the
parabola $(C,C'')$ that we excluded for now, while the second case to the
trajectories of the type $N'NCFO$ starting from a point $N'$ in the region
$\W_4$ and moving with the control $u = -1$ along the parabola $(N',N]$ until
the point $N \in (C, C ''),$ then moving with the control $u = +1$ along the
parabola $[C, C '')$ to the point $C,$ then moving on $M = [t_1, t_2]$
along the state boundary $\G$ to the point $F,$ and finally with the
control $u = +1$ along the curve $S$ to the origin.
In this second case, the value of $\D\mu(t_1)$ can be {\it uniquely determined},
for any normalization of the pair $(\f,\p),$ by the motion time $\D t$ from
the point $N$ to the the point $C$ under the control $u = +1.$
Choosing the normalization $\dot\p = -\f = 1$ on $[0,t_1],$ we have $\p(t_N) =0,$
so $\p(t_C-0) = \D t,$ whence $\D\mu(t_C) = -\D\p(t_C) = \D t.$
\ssk

Thus, for any starting position in $\calD$ not lying on the parabola $(C, C ''),$
there is a unique trajectory satisfying the MP, so it is optimal.\,
\ssk

11) Finally, consider the case of ``boundary'' trajectories, that start from
the points $N$ of parabola $(C, C '')$ on the right of the point $C.$
Any such trajectory can be considered as a tail of a trajectory $N'NCFO$
of the above type.
Since the latter trajectory is optimal, its tail $NCFO$ is also optimal.
\ssk

Note however, that the multipliers in the MP for the trajectory $NCFO,$
unlike those for $N'NCFO,$ are {\it not unique}.\, Indeed, by Lemma \ref{mint2}
we always have $\p=0$ on $[t_1,t_2].$ 
If $\D\mu(t_1) =0,$ then $\p$ and $\f$ are everywhere continuous.
If, moreover, $\f(t_1) <0,$ then $\dot\p= -\f = \const >0$ on $[0,t_1),$
whence $\p<0$ and $u = -1$ there, a contradiction.\, Therefore, $\f(t_1)=0$
and $\f =\p =0$ everywhere on $[0,T],$ a contradiction again.\,

Thus, $\D\mu(t_1)>0,$ and we can take the normalization $\D\mu(t_1) =1.$
Then $\D\p(t_1) = -1,$ hence $\p(t_1-0) =1.$ Herewith $\D\f(t_1) = -q,$ and
since $\f = -\dot\mu \le 0$ on $M,$ then $\f(t_1 +0)= \f(t_1 -0)-q \le 0,$
whence $\f_1 := \f(t_1 -0) \le q.$ Thus, for $t< t_1$ we have
$\f \equiv \f_1 \le q$ and $\p(0) = 1+ \f_1 t_1\,.$ Since $u = +1$ for $t<t_1\,,$
we cannot have $\p<0$ there, so we should have $\p(0)= 1+ \f_1 t_1 \ge 0.$
For $t>t_2$ the measure switches off, $\dot\p = -\f \ge 0,$ so $\p\ge 0,$
which fits to $u =+1$ for any value $\f_1$ bounded only by the inequalities
$-1/t_1 \le \f_1 \le q.$
\ssk

{\it Remark 2.}\, Note the following interesting fact concerning the above
trajectory $NCFO.$ Consider any pair $(\p,\f)$ of costate multipliers
with $\D\mu(t_1) =1,$ satisfying the last inequalities.
If $\f_1 < q,$ then $\f<0$ on $(t_1, t_2]$ and still on $(t_2, T],$ which
yields $\p>0$ on $(t_2, T],$ so $u =+1$ is determined uniquely from the maximum
principle. However, if $\f_1 = q,$ then $\p= 1+ q (t_1-t)$ on $[0,t_1],$
and further $\f= \p =0$ on $(t_1, T].$ In this case $u \in [-1,1]$
{\it can be arbitrary}\, on $(t_1, T],$ with the only requirement that it steers
the point $C$ to the origin.\, Thus, the multipliers with $\f_1=q$ do not select
a unique trajectory satisfying the MP. In particular, these trajectories may
not pass along the boundary interval $(C,F].$ Probably, this peculiar fact
is related to that these trajectories start from the boundary of admissible
set $\calD.$ Obviously, among these trajectories only $NCFO$ is optimal,
because it is a tail of a unique optimal trajectory starting from
a point $N' \in \W_4\,.$ This degeneracy phenomenon deserves more thorough
study in a general setting of problem with state constraints.
\ssk

Thus, the case $k<0$ with three intersection points is completely considered,
and we proved the following

\bth{synth2}
The optimal synthesis is as follows (see Fig. 4). In the sets $\W_1$ and $\W_2$
it is the same as in the free problem, i.e. the optimal trajectory goes
either with the control $u = (-1,+1)$ or with the control $u = (+1,-1).$
From the points of $\W_3$ it moves with control $u = -1$ to the segment $[F, C]$
of the state boundary, then goes along this segment in the ``left--up''
direction until the point $F,$ and then with control $u = +1$ to the origin.
Finally, from the set $\W_4$ the optimal trajectory goes with control $u = -1$
to the parabola $[C, C ''),$ then along this parabola with control $u = +1$ to
the point $C,$ further goes along the segment $[C,F]$ of state boundary, and
at the point $F$ the control switches to $u = +1$ and we get to the origin.
\eth

\subsection{The case $k<0$ with two intersection points}
\label{sec:10}

\begin{figure}[h]
    \begin{center}
        \includegraphics[width=0.6\linewidth,trim={0 0 0 0},clip]{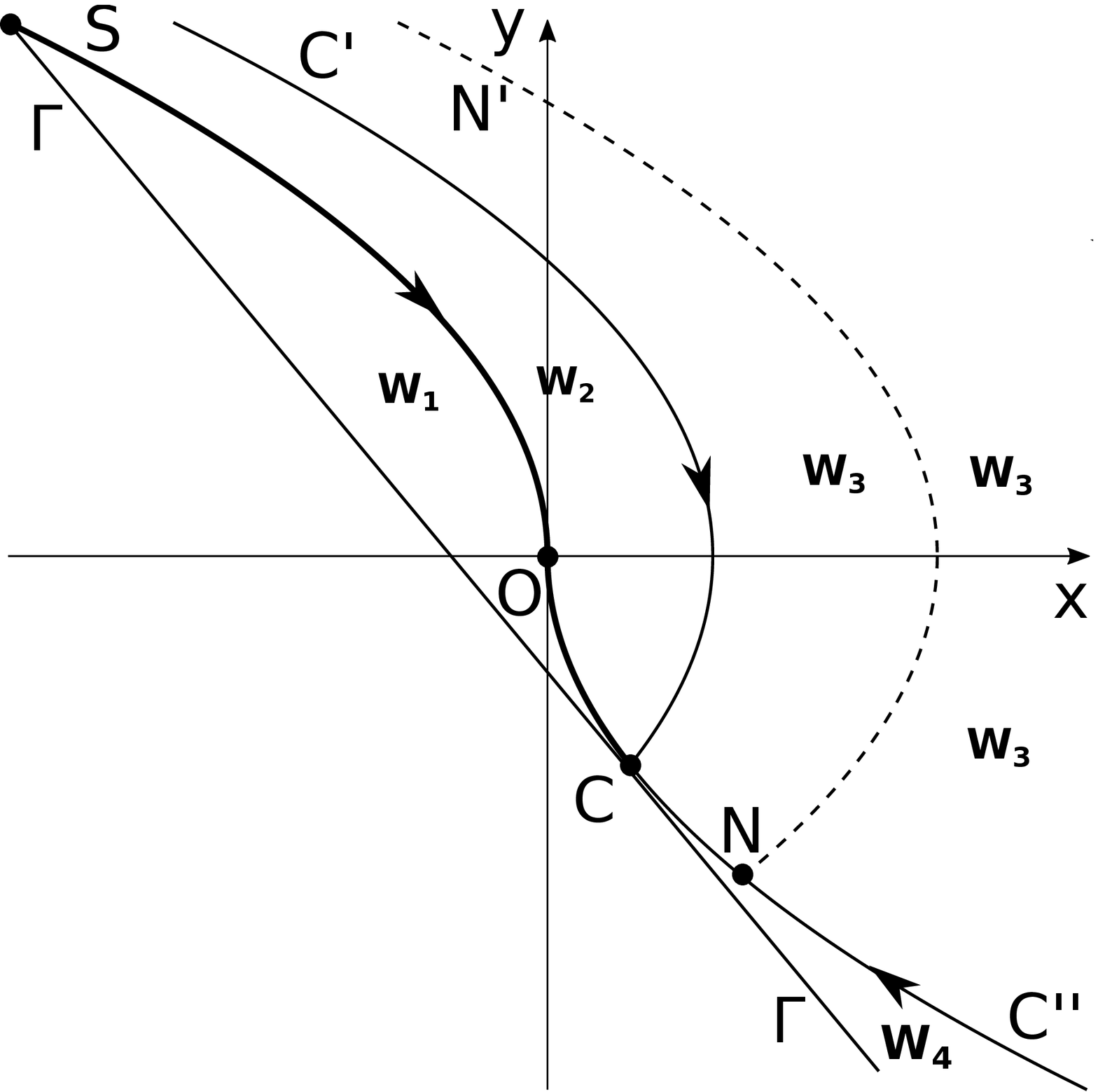}
    \end{center}
    \caption{\label{Fig3}}
\end{figure}

Introduce the characteristic points for this case (see. Fig. \ref{Fig3}).
Let $C$ be the point where $\G$ is tangent to the right branch of the switching
line $S$ (clearly, this point is unique).
\ssk

\noi
Define the following closed sets:

$ W_1$ bounded by the straight line $\G$ and the left branch of the line $S;\,$

$W_2$ bounded by the line $S$ and the parabola $[C, C ');\,$

$W_3$ lying to the right of parabolas $[C, C ')$ and $[C, C'');\,$

and finally,
$W_4$ lying between $\G$ and the parabola $[C, C '').$
\ssk

This case is the limit of the previous case with three intersection points,
when the boundary segment $ [C, E] $ narrows to the single point $C.$
Here the set $W_1$ is the limit of the set $\W_1$ (in some natural sense which
we do not formalize here), the set $W_2$ is the limit of the set
$\W_2 \cap \W_3\,,$ the set $W_3$ is the limit of the set $\W_4\,,$
and the set $W_4$ is the limit of the set $\W_5$.

The points on the state boundary that lie below $C$ are unreachable on
trajectories of our system. \, Any points from $\inter W_4$ are also
unreachable, because any motion from them leads to points of the state boundary
that lie below $C.$ Thus, the set of admissible starting points in this case
is $\calD = W_1 \cup W_2 \cup W_3$.

For the starting points in $\calD,$ the optimal trajectories of the free
problem satisfy everywhere the state constraint, so they remain optimal in
the problem with this constraint as well.
Thus, the optimal synthesis in this case is "the cutting"\, of the free
synthesis along the line composed of the state boundary $\G$ and the
parabola $[C,C '').$
\ssk

{\it Remark 3.}\, Consider the trajectory $C''CO$ starting from a point $C''$
on the parabola $CC''$ and having $u \equiv 1.$ Then $\p\ge 0$ for all $t,$ and
the pair $(\p,\f)$ is not unique. Choose the one with $\p>0$ on $[0,t_C)$ and $\p=0$
on $(t_C, T].$ Setting $\p(0)=1,$ we have $\p(t_C-0) = 1- \f(0)\,t_C \ge 0$
and $\p(t_C+0)=0.$ Since $\D\p(t_C) = -\D\mu(t_C),$ we obtain a unique pair
with $\f(0) = q/(1+q\, t_C).$ Then any admissible trajectory with $u=1$
until $t_C$ and arbitrary after $t_C$ satisfies the MP with these costate
multipliers.

\subsection{The case $k<0$ with one intersection point}
\label{sec:11}

\begin{figure}[h]
    \begin{center}
        \includegraphics[width=0.6\linewidth,trim={0 0 0 0},clip]{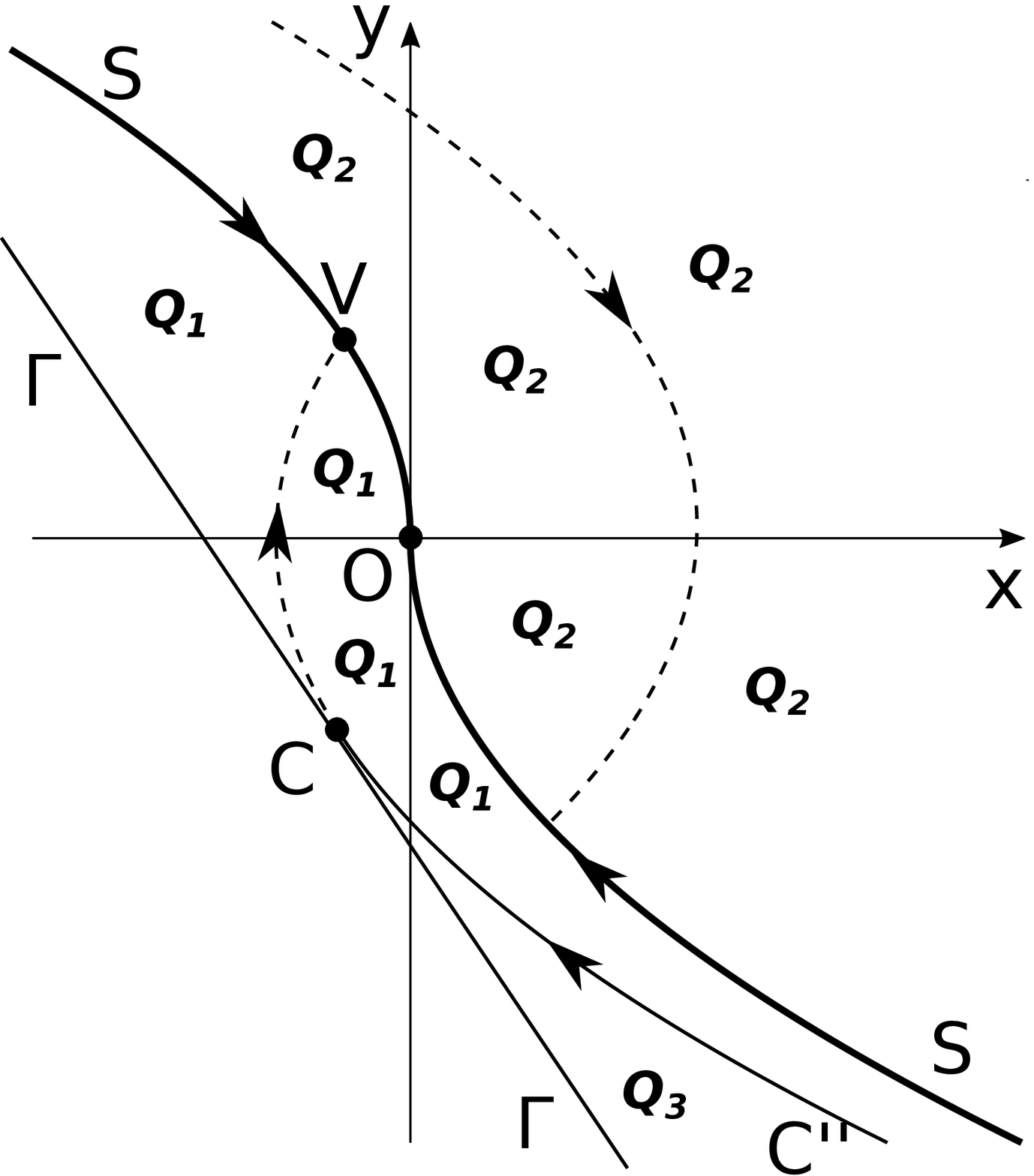}
    \end{center}
    \caption{\label{Fig4}}
\end{figure}

Here, important is the point $C$ where $\G$ is tangent to some parabola
of the family $x= \pol y^2 +\const$ (see Fig. \ref{Fig4}).  Clearly, this
parabola is unique. \ssk

\noi
Define the following closed sets:
\ssk

$Q_1 $ bounded by the line $\G,$ the curve $S,$ and the parabola $[C, C '),\,$

$Q_2$ lying to the right and above the line $\G $ and curve $S,$

and $ Q_3$ lying between the line $\G$ and parabola $[C, C ').$
\ssk

In this case, like before, the points lying on the state boundary below $C,$
as well as points from $\inter Q_3\,,$ are not reachable for the trajectories
of our system. \, Thus, the set of admissible starting points is
$\calD = Q_1 \cup Q_2$.
\ssk

For the starting points in $\calD$ the optimal trajectories of the free
problem satisfy everywhere the state constraint, so they remain optimal in
the problem with this constraint as well.

Thus, the synthesis in this case is "the cutting" \, of the free synthesis
along the line composed by the state boundary $\G $ and the arc of
parabola $[C, C'').$
\ssk

Note that in both three subcases for $k<0,$ the admissible set $\calD$
is closed and the Bellman function is continuous. \ssk

{\it Remark 4.} Consider the trajectory $C''CVO$ starting from a point $C''$
on the parabola $CC''$ and having $u =(1,1,-1).$ Then $\p\ge 0$ for $t<t_V$
and $\p\le0$ for $t>t_V\,,$ while $\D\p(t_C) = -\D\mu \le 0$ and
$\D\dot\p(t_C) = q\,\D\mu \ge 0.$ If $\p(t_C-0)=0,$ then $\D\p(t_C)=0,\,$
$\p=0$ for $t>t_C$ and $\D\dot\p(t_C)=0$ as well, whence $\p=\f=0$ everywhere,
which is forbidden. Therefore, we can set $\p(t_C-0)=1,$ whence
$\p(t_C+0)= 1- \D\mu \ge0.$ Thus, the jump of measure at $t_C$ is arbitrary
within $0\le \D\mu \le 1.$ If $\D\mu= 0,$ then $\p = (t_V-t)/(t_V-t_C)$
is totally linear, and if $\D\mu= 1,$ then $\dot\p= -q$ for $t<t_C\,,$ so
$\p = 1+ q(t_C- t) >0$ there, and $\p=0$ for $t> t_C\,.$

Note that, if $\D\mu < 1,$ the pair $(\p,\f)$ selects the above trajectory
$C''CVO$ only. Otherwise, if $\D\mu =1,$ then $\p>0$ for $t<t_C$ and $\p=\f=0$
for all $t>t_C\,.$ The pair $(\p,\f)$ then selects not the given trajectory only,
but any  admissible one with $u=1$ until $t_C$ and arbitrary $u$ after $t_C\,.$

\section{The Limit Cases}
\label{sec:12}

Consider also the cases when the state boundary $\G$
is horizontal or vertical.

\subsection{The Case $k=0$:\, Horizontal State Boundary }
\label{sec:13}

If $k=0,$ the state constraint is \vad
\begin{equation}\label{state-hor}
y \gee -b \qq (b > 0).
\end{equation}

Here, the admissible set $\calD$ is the whole half-plane $y\ge -b,$
and the analysis of maximum principle is much simpler than that in the
case $k\ne 0.$ By arguments similar to the above ones, it can be easily
shown that the contact set $M$ is still either a singleton or a segment.

\begin{figure}[h!]
    \begin{center}
        \includegraphics[width=0.65\linewidth,trim={0 0 0 0},clip]{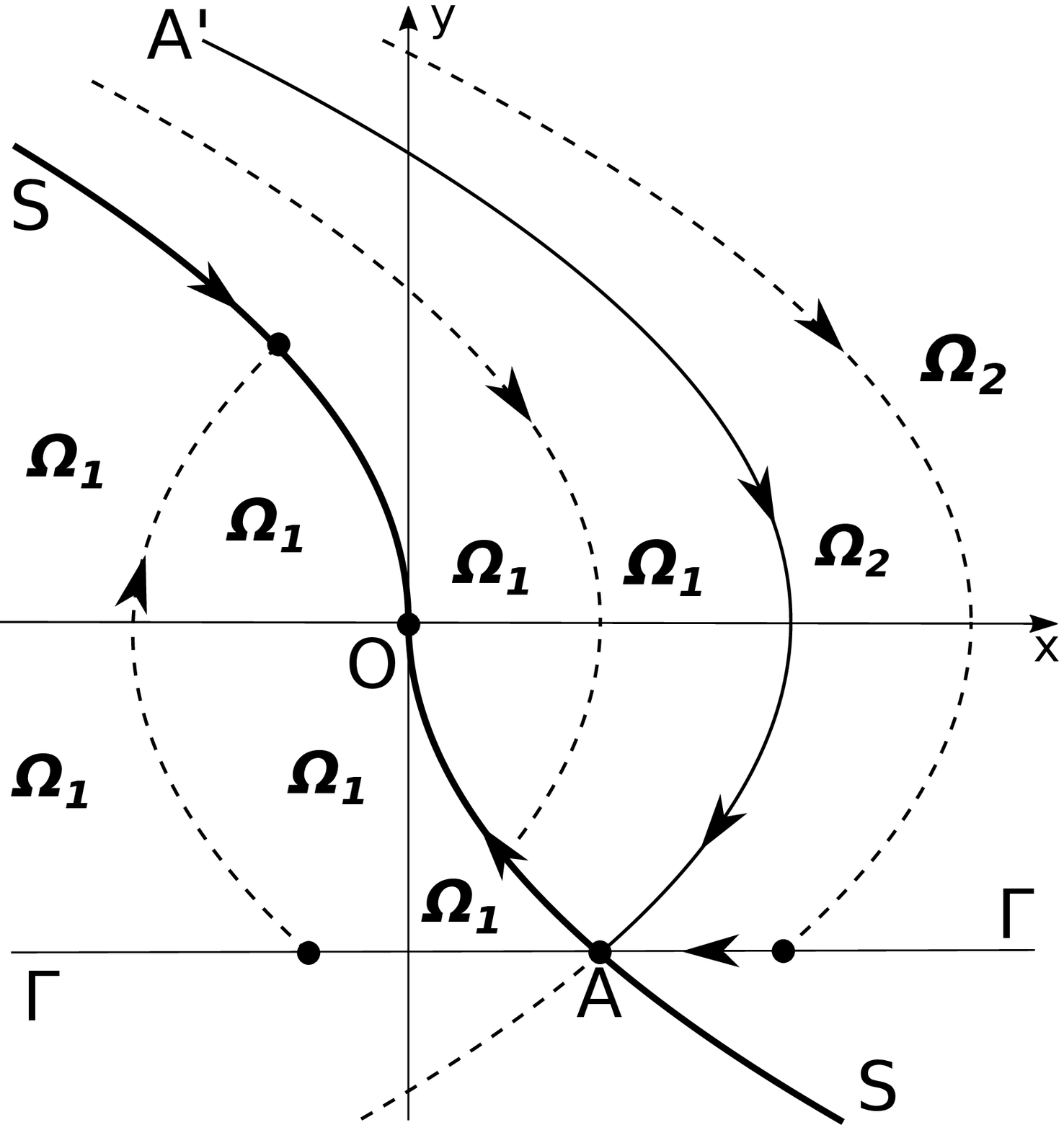}
    \end{center}
    \caption{\label{Fig5}}
\end{figure}

The optimal synthesis is the ``cutting'' of the free synthesis by the state
boundary $\G$ (see Fig. \ref{Fig5}):\, if a trajectory gets to this line from
the domain lying on the right of the switching line $S,$ the control switches
from $u= -1$ to $u=0$ and we move along the line $\G$ up to the point
$A = \G \cap S,$ where the control switches to $u= +1$ and we come to the origin.
\ssk

How this synthesis relates to the synthesis for $k>0,$ shown in
Fig. \ref{Fig1-3}\,?\, If $k\to 0+,$ the point $A$ tends to $(\pol b^2,\, -b),$
the point $C =C(k)$ moves in the left direction (so that $x_C \to -\infty)$
and disappears in the limit. The point $E$ moves to the right along the $x-$axis,
and also disappears in the limit. Therefore, the admissible set $D= D(k)$
transforms in the limit into the half-plane $y\ge -b,$ and the optimal synthesis
for $k>0$ transforms into that for $k=0.$ We do not give here formal definitions
of these limit passages, since their meaning is clear and do not cause
any misunderstandings.
\ssk

Let us consider what happens if $k\to 0-.$ The synthesis for small $k<0$ is shown
in Fig. \ref{Fig2}.\, If $k\to 0-,$ the point $F$ tends to $(\pol b^2,\, -b),$
the tangency point $C =C(k),$ similar to above, moves in the right direction
(so that $x_C \to +\infty)$ and disappears in the limit. Along with this,
the inadmissible domain $\inter\W_5(k)$ moves in the right direction and also
disappears in the limit. One can easily see that here the optimal synthesis
for $k<0$ transforms into that for $k=0.$

\subsection{The Vertical State Boundary $x\le p$}
\label{sec:14}

Consider first the case when the state boundary is vertical and lies on
the right of the origin. This case is a limit one for $k\to +\infty:$
rewriting \eqref{phase} as
$$
\dfrac{y}{k} \gee x - \dfrac{b}{k}\,
$$
and fixing the parameter $p= b/k >0,$ we obtain, in the limit,
the state constraint
\begin{equation}\label{state-vert-1}
x\, \leqslant\, p\,.
\end{equation}

Here, the Pontryagin function and the extended Pontryagin function take
the form, respectively:  \vad
\begin{equation}\label{Pont-vert}
H =\, \varphi y + \psi u, \qq
\overline{H} =\, \varphi y + \psi u - \dot{\mu} \left(x - p\right),
\end{equation}
so that the costate equations now are
\begin{equation}\label{adjEq-vert}
\dot\varphi =\,-\ov H_x =\, \dot\mu, \qq \dot\psi =\,-\ov H_y =\, -\varphi.
\end{equation}

\begin{figure}[h]
    \begin{center}
        \includegraphics[width=0.7\linewidth,trim={0 0 0 0},clip]{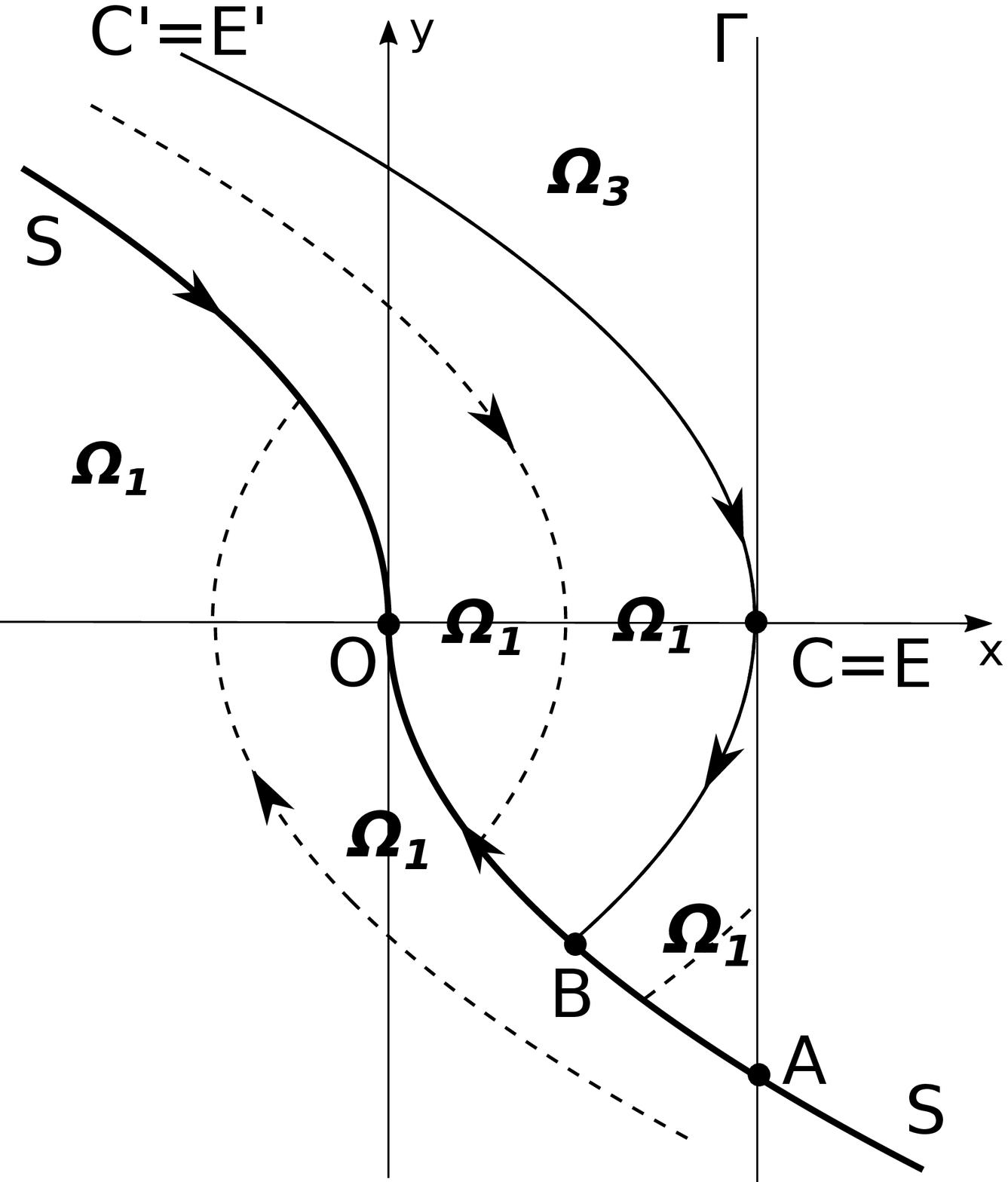}
    \end{center}
    \caption{\label{Fig6}}
\end{figure}

Let us show that here the contact set $M$ with the state boundary can
consist only of the point $C= (p,0),$ i.e. the tangency point of the
state boundary and the parabola $x = p- \pol y^2$ generated by the
control $u= -1$ (see Fig. \ref{Fig6}).
\ssk

Indeed, let at a moment $t'$ we are on $\Gamma$ with $y(t')\neq 0.$
If $y(t') > 0,$ then by the system \eqref{ode} we have $\dx = y>0$
in a neighborhood of $t',$ whence $x>p$ in its right semi-neighborhood,
i.e. we penetrate the state boundary.
This also implies that all the points lying above the parabola $CC'$
and on the left of the line $\G$ are not admissible.

If otherwise $y(t')<0$ for some $t'>0,$ then by the system \eqref{ode}
we have $\dx = y<0$ in a neighborhood of $t',$ whence $x>p$ in its left
semi-neighborhood, so we came to this point from the forbidden zone.

Thus, we can get to the boundary $\G$ for some $t'>0$ only at the point $C,$
and the only way to get to this point is the motion along the tangential
parabola $C'CB$ with control $u= -1.$ In this case, the optimal trajectory
of the free problem, without state constraint, is still admissible in problem
with this constraint, and so, it is still optimal there.

It follows from these considerations that the admissible set $\calD:= \calD(+\infty)$
in this case is bounded from the right by the vertical ray $x=p,\;\, y\le 0$
and the arc of parabola $CC',$ i.e. $\calD(+\infty) =\, \W_1$ in Fig. \ref{Fig6}.
(This figure is the limit of Fig. \ref{Fig1-1} with the set $\W_2$ reduced to
the arc $C'C$ belonging to $\W_1\,.)$
On this set, the optimal synthesis coincides with that for the free problem.

\ssk
However, of interest is the question about possible jump of the measure
$\Delta \mu(t_C)=\delta$ for a trajectory $C'CBO,$ starting from a point
$C'= (x_0, y_0)$ on the arc $C'C$ and having the control $u= (-1, -1, +1).$
Obviously, $y_0 = t_C\,.$

According to the costate equations \eqref{adjEq-vert},
$\p$ is continuous,  \vadjust{\kern-3pt}
$$
\Delta \varphi(t_{C}) =\, \D \mu(t_{C}) =\, \d, \q\, \mbox{and}\q\,
\Delta \dot{\psi}(t_C) =\, - \Delta \f (t_{C})=\, -\delta.
$$
Since $\dot\mu$ can be just an atom at $t_C\,,$ the function $\f$ is piecewise
constant, and $\p$ is piecewise linear with a possible break at the point $t_C\,.$
Obviously, $\p\le 0$ on $[0,t_B]$ and $\p\ge 0$ on $[t_B, T]$ with $\p(t_B)=0.$
If $\p(0)=0,$ then $\p=0$ on the whole interval $[0,T],$ then and also $\f=0$
there, which contradicts the nontriviality.
Therefore, $\p(0)<0,$ so we can set $\psi(0)=\,-1,$ whence $\dot\p >0$ for
$t< t_B\,.$ Since $\D\dot\p(t_C) =\,-\d\le 0$ and $\p(t_B)=0,$
the maximal possible break of $\p$ is when $\p(t_C)=0.$ Then $\dot\p = 1/t_C$
on the left of $t_C$ and $\p= \dot\p =0$ on the right of $t_C,$ so
$\max\,\d = 1/t_C\,.$ Thus, the jump of measure can be arbitrary within
the following bounds:  \vad
$$
0\,\leqslant\, \delta\, \leqslant 1/t_C\, =\, 1/y_0\,.
$$
To verify this, notice that
$\D H(t_C) = \D\f(t_C)\cdot y(t_C) + \p(t_C)\cdot \D u(t_C)=0,$
since $y(t_C) = \p(t_C) =0,$ whence the condition $H= \const$ is valid
for any value of jump $\d.$
\ssk

{\it Remark 5.}
The costate pair $(\p,\f)$ where $\d= 1/t_C\,,$ $\f = -1/t_C$ on $[0,t_C]$
and $\f=0$ on $[t_C, T]\,$ with the corresponding $\p = (t-1)/t_C$ on $[0,t_C]$
and $\p =0$ on $[t_C, T]$ selects not the only trajectory $C'CBO.$
One can easily see that any admissible trajectory passing along the
parabola $C'C$ on $[0,t_C]$ with $u= -1$ and being arbitrary for $t>t_C$
satisfies the maximum principle. However, among these trajectories only the
above $C'CBO$ is optimal, because it is the only optimal one in the free
problem starting from the point $C'.$
\ssk

How the obtained synthesis relates to that in the case $k >> 0,$ the latter
being shown in Fig. \ref{Fig1-1}\,?\, If $k\to +\infty$ with $p = b/k = \const,$
the point $C$ converges in the limit to the point $E = (p, 0).$
Therefore, the set $\W_2(k)$ restricted by the parabolas $CC',$ $EE',$
and the line $\G= \Gamma(k),$ transforms into the arc of parabola $CC',$
and the admissible set $\calD= \calD(k)$ transforms into the set $\calD(+\infty).$
Thus, the optimal synthesis for $k>0$ transforms into the optimal
synthesis for $k= +\infty.$

\subsection{The Vertical State Boundary $\,x\ge -p$}
\label{sec:15}

Consider now the case when $k\to-\infty$ preserving the value
$b/k = -p = \const<0.$ By analogy with the above, passing to the limit,
we get the state constraint  \vad
\begin{equation}\label{state-vert-2}
x\, \ge\, -p \qq (p > 0).
\end{equation}

\begin{figure}[h!]
    \begin{center}
        \includegraphics[width=0.65\linewidth,trim={100 0 0 0},clip]{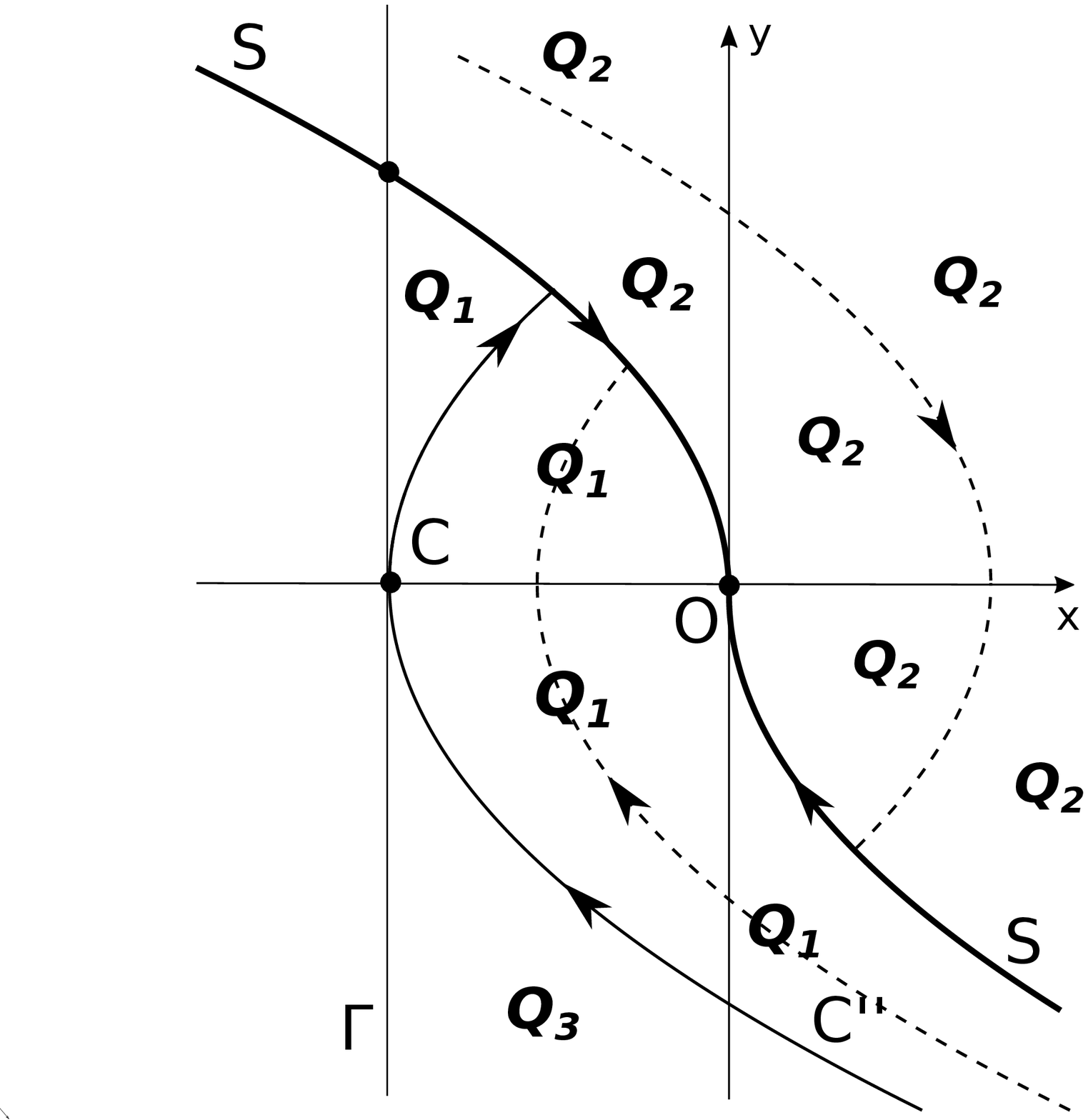}
    \end{center}
    \caption{\label{Fig7}}
\end{figure}

This limit case is symmetrical to the case $x\le p$ with respect to the
origin;\, the corresponding optimal trajectories are shown in Fig. \ref{Fig7}.\,
Here the admissible set $\calD = \calD(-\infty)$ is bounded from the left
by the vertical ray $x= -p,\;\, y\ge0$ and the arc of parabola $CC'',$
i.e. $\calD(-\infty) =\, Q_1 \cup Q_2$ in Fig.~\ref{Fig7}. (The set $Q_3$
is inadmissible.)

Let us consider how this synthesis relates to that for $k<<0,$
shown in Fig. \ref{Fig4}. As $k\to -\infty,$ the touching point $C =C(k)$
converges to $(-p,0),$ the interlocation of the sets $Q_1,\,Q_2,\,Q_3$
remains the same in the limit, and the admissible set $\calD(k)$ transforms
continuously into the above set $\calD(-\infty)$ in Fig.~\ref{Fig7}.
Therefore, the optimal synthesis for $k<< 0$ transforms into that for $k= -\infty.$
\ssk

Note that though the pictures of syntheses for $k= -\infty$ and $k= +\infty$
are symmetric, the prelimit syntheses converge to them in a slightly different
ways:\, in the case $k\to-\infty,$ the structure of admissible starting
positions, shown in Fig. \ref{Fig4}, stays invariable, whereas for $k\to+\infty,$
the region $\W_2(k)$ in Fig. \ref{Fig1-1}, restricted by the parabolas $C C',$
$E E',$ and the line $\G= \Gamma(k),$ disappears in the limit, resulting just
in the arc $C C'$ of tangential parabola in Fig.~\ref{Fig6}.

\section{Conclusion}
We constructed a complete picture of synthesis for the given problem and
analyzed its evolution under the change of parameter $k$ of the slope of
state boundary within the limits $-\infty < k < \infty.\,$
It would be interesting to find the synthesis for the following 3D--problem
with all possible values of $a,b,c$:
$$
\dx=y, \q\; \dy=z,\q\, \dot z=u,\q\, ax +by + cz \ge -d\q (d>0),
$$ $$
|u|\le 1,\q\, (x_0, y_0, z_0)\, \mapsto\, (0,0,0), \q\, T \to \min.
$$

\begin{acknowledgements}
The results for the cases $k>0$ and $k =\pm \infty$ are obtained
by Andrei Dmitruk. His research was supported by the Russian Science Foundation
under grant 20-11-20169 and performed in the Steklov Mathematical Institute
of the Russian Academy of Sciences.\,
The results for the cases $k<0,$ $k\to 0$ and $|k|\to \infty$ are
obtained by Ivan Samylovskiy. His research was supported by the Russian Science
Foundation under grant 19-71-00103 and performed in the Lomonosov Moscow State
University, Faculty of Space Research.\,
The authors thank Valeriy Patsko for useful discussions.
\end{acknowledgements}



\end{document}